\documentclass[oneside,british]{amsart}
\usepackage[T1]{fontenc}
\usepackage[latin9]{inputenc}
\usepackage{amsthm}
\usepackage{amssymb}
\usepackage{setspace}
\usepackage{esint}
\makeatletter
\numberwithin{equation}{section}
\numberwithin{figure}{section}
  \theoremstyle{plain}
  \newtheorem*{thm*}{\protect\theoremname}
  \theoremstyle{remark}
  \newtheorem*{rem*}{\protect\remarkname}
\theoremstyle{plain}
\newtheorem{thm}{\protect\theoremname}
  \theoremstyle{definition}
  \newtheorem{defn}[thm]{\protect\definitionname}
  \theoremstyle{remark}
  \newtheorem{rem}[thm]{\protect\remarkname}
  \theoremstyle{plain}
  \newtheorem{prop}[thm]{\protect\propositionname}
  \theoremstyle{plain}
  \newtheorem{lem}[thm]{\protect\lemmaname}

\date{}

\makeatother

\usepackage{babel}
  \providecommand{\definitionname}{Definition}
  \providecommand{\lemmaname}{Lemma}
  \providecommand{\propositionname}{Proposition}
  \providecommand{\remarkname}{Remark}
  \providecommand{\theoremname}{Theorem}
\providecommand{\theoremname}{Theorem}

\begin{document}

\title{On self-affine measures with equal Hausdorff and Lyapunov dimensions}

\author{Ariel Rapaport}

\date{October 24, 2015}



\thanks{Supported by ERC grant 306494}
\begin{abstract}
Let $\mu$ be a self-affine measure on $\mathbb{R}^{d}$ associated
to a self-affine IFS $\{\varphi_{\lambda}(x)=A_{\lambda}x+v_{\lambda}\}_{\lambda\in\Lambda}$
and a probability vector $p=(p_{\lambda})_{\lambda}>0$. Assume the
strong separation condition holds. Let $\gamma_{1}\ge...\ge\gamma_{d}$
and $D$ be the Lyapunov exponents and dimension corresponding to
$\{A_{\lambda}\}_{\lambda\in\Lambda}$ and $p^{\mathbb{N}}$, and
let $\mathbf{G}$ be the group generated by $\{A_{\lambda}\}_{\lambda\in\Lambda}$.
We show that if $\gamma_{m+1}>\gamma_{m}=...=\gamma_{d}$, if $\mathbf{G}$
acts irreducibly on the vector space of alternating $m$-forms, and
if the Furstenberg measure $\mu_{F}$ satisfies $\dim_{H}\mu_{F}+D>(m+1)(d-m)$,
then $\mu$ is exact dimensional with $\dim\mu=D$.
\end{abstract}

\maketitle

\section{Introduction}

Let $d\ge2$ and let $\Lambda$ be a finite index set. Fix a family
of matrices $\{A_{\lambda}\}_{\lambda\in\Lambda}=\mathbf{A}\subset Gl(d,\mathbb{R})$
with $\left\Vert A_{\lambda}\right\Vert <1$ for $\lambda\in\Lambda$,
let $\{v_{\lambda}\}_{\lambda\in\Lambda}\subset\mathbb{R}^{d}$, and
fix a probability vector $p=\{p_{\lambda}\}_{\lambda\in\Lambda}>0$.
Let $\{\varphi_{\lambda}\}_{\lambda\in\Lambda}$ be the self-affine
IFS with
\begin{equation}
\varphi_{\lambda}(x)=A_{\lambda}x+v_{\lambda}\:\mbox{ for }\lambda\in\Lambda\mbox{ and }x\in\mathbb{R}^{d}\:.\label{E25}
\end{equation}
Denote by $\mu$ the self-affine measure on $\mathbb{R}^{d}$ which
corresponds to $\{\varphi_{\lambda}\}_{\lambda\in\Lambda}$ and $p$,
i.e. $\mu$ is the unique probability measure with
\[
\mu=\sum_{\lambda\in\Lambda}p_{\lambda}\cdot\varphi_{\lambda}\mu\:.
\]
The Lyapunov dimension $D$ of $\mu$ (see Section \ref{S2} below)
is an upper bound for the dimension of $\mu$, but it is in general
difficult to verify whether there is equality. The purpose of this
paper is to present verifiable conditions under which
\begin{equation}
\mu\mbox{ is exact dimensional with }\dim\mu=D\:.\label{E22}
\end{equation}

\subsection{Background for the problem}

Let us mention some notable results regarding self-affine measures
and sets. From Theorem 1.9 in \cite{JPS} it follows that $D$ is
the 'typical' value of $\dim_{H}\mu$, where $\dim_{H}$ stands for
the Hausdorff dimension. More precisely, it is shown that if $\left\Vert A_{\lambda}\right\Vert <\frac{1}{2}$
for $\lambda\in\Lambda$ and if the translations $\{v_{\lambda}\}_{\lambda\in\Lambda}$
are drawn according to the Lebesgue measure, then $\dim_{H}\mu=\min\{D,d\}$
almost surely. The inequality $\dim_{H}\mu\le D$ is always satisfied.

Analogous to this is the following classical result, due to Falconer,
regarding the typical dimension of self-affine sets. Let $K$ be the
attractor of $\{\varphi_{\lambda}\}_{\lambda\in\Lambda}$. In \cite{F3}
it is shown that if $\left\Vert A_{\lambda}\right\Vert <\frac{1}{3}$
for $\lambda\in\Lambda$, then
\[
\dim_{H}K=\min\{\dim_{A}K,d\}\mbox{ for Lebesgue almost all }\{v_{\lambda}\}_{\lambda\in\Lambda}\:.
\]
Here $\dim_{A}K$ stands for the affinity dimension of $K$, which
is defined in terms of the matrices in $\mathbf{A}$. This was later
improved in \cite{S} by replacing the constant $\frac{1}{3}$ by
$\frac{1}{2}$. The inequality $\dim_{H}K\le\dim_{A}K$ is always
true.

For fixed translations $\{v_{\lambda}\}_{\lambda\in\Lambda}$ the
exact value of $\dim_{H}K$ has been found for several specific classes
of self-affine sets. See the survey \cite{F4} and the references
therein. Much attention has been given to fractal carpets, where members
of $\mathbf{A}$ preserve horizontal and vertical directions (see
\cite{M1} for instance).

Here we establish (\ref{E22}) in the opposite situation, in which
there is no proper subspace invariant under all members of $\mathbf{A}$.
This makes it possible to consider the Furstenberg measure $\mu_{F}$
on the Grassmannian manifold (see Section \ref{S2} below). The measure
$\mu_{F}$ allows us to control the distribution of the orientation
of cylinder sets at small scale.

For $d=2$ this idea was already used in \cite{FK} and \cite{B},
in order to obtain (\ref{E22}) under assumptions different than ours.
In Section \ref{S1.4} below we describe these results and compare
them with the work presented here. A notable advantage in our result
is that we do not require a lower bound on $\dim_{H}\mu$, but rather
only on $D$ which is at least as large and independent of the translations
$\{v_{\lambda}\}_{\lambda\in\Lambda}$.

\subsection{The main result}

We shall consider only the case where the IFS $\{\varphi_{\lambda}\}_{\lambda\in\Lambda}$
satisfies the strong separation condition (SSC). Denote by $\gamma_{1}\ge...\ge\gamma_{d}$
the Lyapunov exponents corresponding to the Bernoulli measure $p^{\mathbb{N}}$
and the matrices $\mathbf{A}$, and set
\[
m=\max\{1\le i\le d\::\:\gamma_{d-i+1}=...=\gamma_{d}\}\:.
\]
If $m=d$ and the SSC is satisfied then (\ref{E22}) follows directly
from Theorem 2.6 in \cite{FH}. Hence assume $m<d$. Let $\mathbf{G}\subset Gl(d,\mathbb{R})$
be the closure of the group generated by $\mathbf{A}$. We assume
that $\mathbf{G}$ is $m$-irreducible, which means that it acts irreducibly
on the vector space of alternating $m$-forms. A precise definition
is given in Section \ref{S2}. When $m=1\mbox{ or }d-1$, and in particular
when $d=2\mbox{ or }3$, this condition reduces to the absence of
a proper subspace of $\mathbb{R}^{d}$ which is invariant under all
members of $\mathbf{A}$ (see remark \ref{R2} below).

Let $G_{d,m}$ denote the Grassmannian manifold of all $m$-dimensional
linear subspaces of $\mathbb{R}^{d}$. Each $M\in Gl(d,\mathbb{R})$
defines a map from $G_{d,m}$ onto itself, which takes $W\in G_{d,m}$
to $M(W)$. From $m<d$, the irreducibility assumption, and results
found in \cite{BL}, it follows that there exists a unique probability
measure $\mu_{F}$ on $G_{d,m}$ with
\[
\mu_{F}=\sum_{\lambda\in\Lambda}p_{\lambda}\cdot A_{\lambda}^{-1}\mu_{F},
\]
and moreover that $\dim_{H}\mu_{F}>0$ (see Proposition \ref{L2}
in Section \ref{S2}). The measure $\mu_{F}$ is called the Furstenberg
measure on $G_{d,m}$ corresponding to $\mathbf{A}^{-1}:=\{A_{\lambda}^{-1}\}_{\lambda\in\Lambda}$
and $p$. The following theorem is our main result.
\begin{thm*}
Assume the following conditions:$\newline$(i) $\{\varphi_{\lambda}\}_{\lambda\in\Lambda}$
satisfies the SSC,$\newline$(ii) $m$ is strictly smaller than $d$,$\newline$(iii)
$\mathbf{G}$ is $m$-irreducible, and$\newline$(iv) The measure
$\mu_{F}$ satisfies
\[
\dim_{H}\mu_{F}+D>(m+1)(d-m)\:.
\]
Then (\ref{E22}) holds true, i.e. $\mu$ is exact dimensional with
$\dim\mu=D$.
\end{thm*}

\subsection{Explicit examples}

The theorem just stated can be used to compute the dimension of many
concrete self-affine measures. In order to do so one needs to bound
$\dim_{H}\mu_{F}$ from below, which is not a trivial problem. Let
us mention some results which are relevant for this task. Here we
assume the elements of $\mathbf{A}$ are distinct, i.e. $A_{\lambda_{1}}\ne A_{\lambda_{2}}$
for $\lambda_{1},\lambda_{2}\in\Lambda$ with $\lambda_{1}\ne\lambda_{2}$.
Also, we shall have no need for the matrices in $\mathbf{A}$ to be
contractions. Indeed, the Furstenberg measure is unaffected if we
multiply members of $\mathbf{A}$ by non-zero scalars.

In \cite{HS} it is shown that if $\mathbf{A}\subset Gl(2,\mathbb{R})$
and $p$ are such that elements in $\mathbf{A}$ have algebraic entries
and determinant $1$, $\mathbf{A}$ generates a free group, $\gamma_{1}$
is strictly greater than $\gamma_{2}$, and $\mathbf{G}$ acts irreducibly
on $\mathbb{R}^{2}$, then 
\[
\dim_{H}\mu_{F}=\min\{\frac{H(p)}{-2\cdot\gamma_{1}},1\}\:.
\]
Here $H(p)$ stands for the entropy of $p$. For example, this can
be applied when $p>0$ and
\begin{equation}
\mathbf{A}=\left\{ \left(\begin{array}{cc}
1 & 2\\
0 & 1
\end{array}\right),\left(\begin{array}{cc}
1 & 0\\
2 & 1
\end{array}\right)\right\} \:.\label{E23}
\end{equation}

In Section VI.5 of \cite{BL} it is shown that $\dim_{H}\mu_{F}=\frac{H(p)}{-2\cdot\gamma_{1}}$
whenever $|\mathbf{A}|>1$, $p>0$, and
\[
\mathbf{A}^{-1}\subset\left\{ \left(\begin{array}{cc}
0 & 1\\
1 & n
\end{array}\right)\::\:n\ge1\right\} \:.
\]

For $E,L\in\mathbb{R}$ with $|E|+|L|<2$, denote by $\mu_{F}^{E,L}$
the Furstenberg measure corresponding to
\begin{equation}
\mathbf{A}^{-1}=\left\{ \left(\begin{array}{cc}
E-L & -1\\
1 & 0
\end{array}\right),\left(\begin{array}{cc}
E+L & -1\\
1 & 0
\end{array}\right)\right\} \mbox{ and }p=(\frac{1}{2},\frac{1}{2})\:.\label{E24}
\end{equation}
In \cite{B2} it is shown that there exists a constant $\delta>0$
with
\[
\underset{L\rightarrow0}{\lim}\:\dim_{H}\mu_{F}^{E,L}=1\:\mbox{ for all }\:\delta<|E|<2-\delta\:.
\]

In \cite{B3} an example is given, for the case $d=2$, of $\mathbf{A}$
and $p$ for which $\gamma_{1}>\gamma_{2}$, the action of $\mathbf{G}$
is irreducible, and $\mu_{F}$ is absolutely continuous with respect
to the Lebesgue measure. For $d\ge3$ an example of $\mathbf{A}$
and $p$ with these properties is obtained in \cite{BQ2}.

\subsection{\label{S1.4}Comparison with recent work}

As mentioned above, for $d=2$ the validity of (\ref{E22}) was established
in two recent papers under conditions different than ours. From the
arguments found in \cite{FK}, it follows that if the matrices in
$\mathbf{A}$ have strictly positive entries, $\{\varphi_{\lambda}\}_{\lambda\in\Lambda}$
satisfies the SSC, and 
\[
\dim_{H}\mu_{F}+\dim_{H}\mu>2,
\]
then (\ref{E22}) holds. This is actually done more generally, in
the sense that the self-affine measure $\mu$ can be replaced by the
projection of a Gibbs measures into $\mathbb{R}^{2}$.

Given $M\in Gl(2,\mathbb{R})$ let $\alpha_{1}(M)\ge\alpha_{2}(M)>0$
denote the singular values of $M$. It is said that $\mathbf{A}$
satisfies the dominated splitting condition if there exist constants
$0<C,\delta<\infty$ with
\[
\frac{\alpha_{1}(A_{1}\cdot...\cdot A_{n})}{\alpha_{2}(A_{1}\cdot...\cdot A_{n})}\ge C\cdot e^{\delta n}\mbox{ for all }n\ge1\mbox{ and }A_{1},...,A_{n}\in\mathbf{A}\:.
\]
For example, this is satisfied when the matrices in $\mathbf{A}$
have strictly positive entries. It is shown in \cite{B} that if $\mathbf{A}$
satisfies dominated splitting, $\{\varphi_{\lambda}\}_{\lambda\in\Lambda}$
satisfies the SSC, and 
\[
\dim_{H}\mu_{F}+\dim_{H}\mu>2\:\mbox{ or }\:\dim_{H}\mu_{F}\ge\min\{1,D\},
\]
then (\ref{E22}) holds.

Note that since $D\ge\dim_{H}\mu$, the condition $\dim_{H}\mu_{F}+D>2$,
which appears in our result when $d=2$, is weaker than $\dim_{H}\mu_{F}+\dim_{H}\mu>2$.
This is important because $D$, as opposed to $\dim_{H}\mu$, is independent
of the choice of translations $\{v_{\lambda}\}_{\lambda\in\Lambda}$.
Observe also, that if the closure of the set
\[
\{A_{1}\cdot...\cdot A_{n}\::\:n\ge1\mbox{ and }A_{1},...,A_{n}\in\mathbf{A}\}
\]
contains an element $A\in Gl(2,\mathbb{R})$ for which $\frac{\alpha_{1}(A^{n})}{\alpha_{2}(A^{n})}$
does not increase exponentially fast as $n\rightarrow\infty$, then
the results from \cite{B} and \cite{FK} don't apply but our result
can. This is in fact the case in examples (\ref{E23}) and (\ref{E24})
mentioned above. This is also true for the example obtained in \cite{B3},
since in that case $A^{-1}\in\mathbf{A}$ whenever $A\in\mathbf{A}$.

By using the aforementioned results about measures, results about
the dimension of certain self-affine sets are obtained in \cite{B}
and \cite{FK}. More precisely, conditions for $\dim_{H}K=\dim_{A}K$
are given, where recall that $K$ is the attractor of $\{\varphi_{\lambda}\}_{\lambda\in\Lambda}$
and $\dim_{A}K$ is the affinity dimension of $K$. We do not pursue
this here, although it seems reasonable to believe that our work can
also be applied in order to obtain this equality for new classes of
self-affine sets.
\begin{rem*}
In the last stages of writing up this research the author became aware
of the preprint \cite{BK}. When $d=2$ it is shown in \cite{BK}
that $\mu$ is always exact dimensional, and that $\dim\mu=D$ if
the SSC holds and 
\[
\dim_{H}\mu_{F}>\min\{\dim\mu,\:2-\dim\mu\}\:.
\]
As mentioned above, since $D\ge\dim_{H}\mu$ our result may be easier
to use in some cases. For $d>2$ results are proven in \cite{BK}
under an assumption on $\mathbf{A}$, termed totally dominated splitting,
which is a multi-dimensional analogue of the dominated splitting condition
previously mentioned. Hence for $d>2$ our work applies in many situations
that are untreated by \cite{BK}.
\end{rem*}

\subsection{About the proof}

We now make the dependency in the translations explicit. Given $(v_{\lambda})_{\lambda\in\Lambda}=v\in\mathbb{R}^{d|\Lambda|}$
denote by $\{\varphi_{v,\lambda}\}_{\lambda\in\Lambda}$ the IFS satisfying
(\ref{E25}), and let $\mu_{v}$ be the self-affine measure corresponding
to $\{\varphi_{v,\lambda}\}_{\lambda\in\Lambda}$ and $p$. Let $\mathcal{V}\subset\mathbb{R}^{d|\Lambda|}$
be the set of all $v\in\mathbb{R}^{d|\Lambda|}$ for which $\{\varphi_{v,\lambda}\}_{\lambda\in\Lambda}$
satisfies the SSC. In the proofs found in \cite{B} and \cite{FK},
some $v\in\mathcal{V}$ is fixed and linear projections and sections
of the measure $\mu_{v}$ are studied. In our proof we shall also
examine linear sections of measures, but we shall consider the entire
collection $\{\mu_{v}\}_{v\in\mathcal{V}}$ at once.

More precisely, it will be shown that there exists an upper semi-continuous
function $F:\mathcal{V}\rightarrow[0,\infty)$, such that for every
$v\in\mathcal{V}$ and $\mu_{v}\times\mu_{F}$-a.e. $(x,W)\in\mathbb{R}^{d}\times G_{d,m}$
the sliced measure, obtained from $\mu_{v}$ and supported on $x+W$,
has exact dimension $F(v)$. The proof of this uses ergodic theory
and results from the random matrix theory presented in \cite{BL}.
From the result of \cite{JPS} mentioned above, and from results found
in \cite{M2} regarding the dimension of exceptional sets of sections,
it will follows that $F(v)\ge D-d+m$ for $\mathcal{L}eb$-a.e. $v\in\mathcal{V}$.
The semi-continuity of $F$ implies that this inequality holds in
fact for all $v\in\mathcal{V}$. Now by fixing $v\in\mathcal{V}$
and using estimates on the dimension of exceptional sets of projections,
it will follows that $\dim_{H}\mu_{v}\ge D$. The inequality $\dim\mu_{v}\le D$
in not hard to prove, and completes the proof.

\subsection{Outline of the paper}

In Section \ref{S2} we give some necessary definitions and state
Theorem \ref{T3} which is our main result. Is Section \ref{S3} we
carry out the proof, while relaying on Proposition \ref{P3} and Lemmas
\ref{L6} to \ref{L11}, whose proofs are deferred to subsequent sections.
In Section \ref{S4} we state and prove some required results, which
follow from the theory of random matrices. In Section \ref{S5} we
prove Proposition \ref{P3}, which is the main ingredient in the proof
of Theorem \ref{T3}. In Section \ref{S6} we prove all auxiliary
lemmas which were priorly used without proof.

\textbf{Acknowledgement. }I would like to thank my advisor Michael
Hochman, for suggesting to me the problem studied in this paper, and
for many helpful discussions.

\section{\label{S2}Statement of the main result}

Fix some integer $d\ge2$ and for $x\in\mathbb{R}^{d}$ denote by
$|x|$ the euclidean norm of $x$. For a $d\times d$ matrix $M$
(or operator on $\mathbb{R}^{d}$) denote by $\left\Vert M\right\Vert $
the operator norm of $M$ with respect to the euclidean norm. Let
$\Lambda$ be a finite set with $|\Lambda|>1$, and fix $\{A_{\lambda}\}_{\lambda\in\Lambda}\subset Gl(d,\mathbb{R})$
with $\left\Vert A_{\lambda}\right\Vert <1$ for each $\lambda\in\Lambda$.
Let $\mathbf{G}\subset Gl(d,\mathbb{R})$ be the closure of the group
generated by $\{A_{\lambda}\}_{\lambda\in\Lambda}$. For $(v_{\lambda})_{\lambda\in\Lambda}=v\in\mathbb{R}^{d|\Lambda|}$
let $\{\varphi_{v,\lambda}\}_{\lambda\in\Lambda}$ be the self-affine
IFS with $\varphi_{v,\lambda}(x)=A_{\lambda}x+v_{\lambda}$ for $\lambda\in\Lambda$
and $x\in\mathbb{R}^{d}$. Let $K_{v}\subset\mathbb{R}^{d}$ be the
attractor of $\{\varphi_{v,\lambda}\}_{\lambda\in\Lambda}$, i.e.
$K_{v}$ is the unique non empty compact subset of $\mathbb{R}^{d}$
with $K_{v}=\cup_{\lambda\in\Lambda}\varphi_{v,\lambda}(K_{v})$.
We say that the strong separation condition (SSC) holds for $\{\varphi_{v,\lambda}\}_{\lambda\in\Lambda}$
if the union $\cup_{\lambda\in\Lambda}\varphi_{v,\lambda}(K_{v})$
is disjoint, and we denote by $\mathcal{V}\subset\mathbb{R}^{d|\Lambda|}$
the set of all $v\in\mathbb{R}^{d|\Lambda|}$ for which the SSC holds.
It is easy to see that $\mathcal{V}$ is an open subset of $\mathbb{R}^{d|\Lambda|}$,
and we assume it to be non empty.

Let $p=(p_{\lambda})_{\lambda\in\Lambda}$ be a probability vector
with $p_{\lambda}>0$ for each $\lambda\in\Lambda$. Set $\Omega=\Lambda^{\mathbb{N}}$,
equip $\Lambda$ with the discrete topology, and equip $\Omega$ with
the product topology. Let $\mathcal{F}$ be the Borel $\sigma$-algebra
of $\Omega$, and let $\mu$ be the Bernoulli measure on $(\Omega,\mathcal{F})$
which corresponds to $p$ (i.e. $\mu=p^{\mathbb{N}}$). For each $v\in\mathbb{R}^{d|\Lambda|}$
and $\omega\in\Omega$ set 
\[
\pi_{v}(\omega)=\underset{n}{\lim}\:\varphi_{v,\omega_{0}}\circ...\circ\varphi_{v,\omega_{n}}(0)\:.
\]
Since the mappings $\{\varphi_{v,\lambda}\}_{\lambda\in\Lambda}$
are contractions this limit always exists and $\pi_{v}:\Omega\rightarrow\mathbb{R}^{d}$
is continuous. Note that $\pi_{v}\mu:=\mu\circ\pi_{v}^{-1}$ is the
unique Borel probability measure on $\mathbb{R}^{d}$ for which the
relation $\pi_{v}\mu=\sum_{\lambda\in\Lambda}p_{\lambda}\cdot\varphi_{v,\lambda}\pi_{v}\mu$
is satisfied.

Given $M\in Gl(d,\mathbb{R})$ let $\alpha_{1}(M)\ge...\ge\alpha_{d}(M)>0$
be the singular values of $M$. Let $0>\gamma_{1}\ge...\ge\gamma_{d}>-\infty$
be the Lyapunov exponents corresponding to $\mu$ and $\{A_{\lambda}\}_{\lambda\in\Lambda}$
(see chapter III.5 in \cite{BL}), i.e. for $\mu$-a.e. $\omega\in\Omega$
\begin{equation}
\gamma_{i}=\underset{n}{\lim}\frac{1}{n}\log\alpha_{i}(A_{\omega_{0}}\cdot...\cdot A_{\omega_{n-1}})\mbox{ for }1\le i\le d\:.\label{E1}
\end{equation}
Denote the entropy of $\mu$ by $h_{\mu}$ (i.e. $h_{\mu}=\sum_{\lambda\in\Lambda}-p_{\lambda}\cdot\log p_{\lambda}$),
set
\begin{equation}
k(\mu)=\max\{0\le i\le d\::\:0<h_{\mu}+\gamma_{1}+...+\gamma_{i}\},\label{E19}
\end{equation}
and set
\[
D(\mu)=\begin{cases}
k(\mu)-\frac{h_{\mu}+\gamma_{1}+...+\gamma_{k(\mu)}}{\gamma_{k(\mu)+1}} & ,\mbox{ if }k(\mu)<d\\
-d\cdot\frac{h_{\mu}}{\gamma_{1}+...+\gamma_{d}} & ,\mbox{ }\mbox{if }k(\mu)=d
\end{cases}\:.
\]
The number $D(\mu)$ is called the Lyapunov dimension of $\mu$ with
respect to the family $\{A_{\lambda}\}_{\lambda\in\Lambda}$.

Given a metric space $X$ we denote the collection of all compactly
supported Borel probability measures on $X$ by $\mathcal{M}(X)$.
For $\theta\in\mathcal{M}(X)$ we write
\[
\dim_{H}\theta=\inf\{\dim_{H}E\::\:E\subset X\mbox{ is a Borel set with }\theta(E)>0\}
\]
and
\[
\dim_{H}^{*}\theta=\inf\{\dim_{H}E\::\:E\subset X\mbox{ is a Borel set with }\theta(X\setminus E)=0\},
\]
where $\dim_{H}E$ stands for the Hausdorff dimension of the set $E$.
For $x\in\mathbb{R}^{d}$ and $\epsilon>0$ denote by $B(x,\epsilon)$
the closed ball in $\mathbb{R}^{d}$ with centre $x$ and radius $\epsilon$.
Given $\theta\in\mathcal{M}(\mathbb{R}^{d})$ we say that $\theta$
has exact dimension $s\ge0$ if
\[
\underset{\epsilon\downarrow0}{\lim}\frac{\log\theta(B(x,\epsilon))}{\log\epsilon}=s\mbox{ for }\theta\mbox{-a.e. }x\in\mathbb{R}^{d},
\]
in which case we write $\dim\theta=s$. It is well known (see chapter
10 of \cite{F}) that
\begin{equation}
\dim_{H}\theta=\mathrm{essinf}_{\theta}\{\underset{\epsilon\downarrow0}{\liminf}\frac{\log\theta(B(x,\epsilon))}{\log\epsilon}\::\:x\in\mathbb{R}^{d}\}\:.\label{E2}
\end{equation}

Given $1\le m<d$ let $G_{d,m}$ denote the Grassmannian manifold
of all $m$-dimensional linear subspaces of $\mathbb{R}^{d}$. For
a subspace $W\subset\mathbb{R}^{d}$ let $P_{W}:\mathbb{R}^{d}\rightarrow\mathbb{R}^{d}$
be the orthogonal projection onto $W$. For $W,U\in G_{d,m}$ set
$d_{G_{d,m}}(W,U)=\left\Vert P_{W}-P_{U}\right\Vert $, then $d_{G_{d,m}}$
is a metric on $G_{d,m}$ which we shall use. For $M\in Gl(d,\mathbb{R})$
and $W\in G_{d,m}$ set $M\cdot W=M(W)\in G_{d,m}$, which defines
an action of $Gl(d,\mathbb{R})$ on $G_{d,m}$.

For $1\le m\le d$ let $\mathcal{A}^{m}(\mathbb{R}^{d})$ denote the
vector space of alternating $m$-linear forms on $(\mathbb{R}^{d})^{*}$.
Given $x_{1},...,x_{m}\in\mathbb{R}^{d}$ let $x_{1}\wedge...\wedge x_{m}\in\mathcal{A}^{m}(\mathbb{R}^{d})$
be such that 
\[
x_{1}\wedge...\wedge x_{m}(f_{1},...,f_{m})=\det[\{f_{i}(x_{j})\}_{i,j=1}^{m}]\mbox{ for }f_{1},...,f_{m}\in(\mathbb{R}^{d})^{*}\:.
\]
If $\{e_{1},...,e_{n}\}$ is a basis for $\mathbb{R}^{d}$ then
\[
\{e_{i_{1}}\wedge...\wedge e_{i_{m}}\::\:1\le i_{1}<...<i_{m}\le d\}
\]
is a basis for $\mathcal{A}^{m}(\mathbb{R}^{d})$. For $M\in Gl(d,\mathbb{R})$
we define an automorphism $\mathcal{A}^{m}M$ of $\mathcal{A}^{m}(\mathbb{R}^{d})$
by
\[
\mathcal{A}^{m}M(x_{1}\wedge...\wedge x_{m})=Mx_{1}\wedge...\wedge Mx_{m}\mbox{ for }x_{1},...,x_{m}\in\mathbb{R}^{d}\:.
\]

\begin{defn}
Given $1\le m<d$ and $\mathbf{S}\subset Gl(d,\mathbb{R})$ we say
that $\mathbf{S}$ is $m$-irreducible if there does not exist a proper
linear subspace $W$ of $\mathcal{A}^{m}(\mathbb{R}^{d})$ with $\mathcal{A}^{m}M(W)=W$
for each $M\in\mathbf{S}$. When $m=1$ we say that $\mathbf{S}$
is irreducible.
\end{defn}

\begin{rem}
\label{R2}Clearly $\mathbf{S}$ is irreducible if and only if there
does not exist a proper linear subspace $W$ of $\mathbb{R}^{d}$
with $M(W)=W$ for each $M\in\mathbf{S}$. It is also easy to show
that $\mathbf{S}$ is $m$-irreducible if and only if it is $d-m$-irreducible
(see page 86 in \cite{BL}). Hence when $d=2\mbox{ or }3$ the $m$-irreducibility
condition reduces to the absence of a proper subspace of $\mathbb{R}^{d}$
which is $M$-invariant for all $M\in\mathbf{S}$.
\end{rem}

The following proposition follows from results found in \cite{BL},
and shall be proven in Section \ref{S4}. From now on we set
\[
m=\max\{1\le i\le d\::\:\gamma_{d-i+1}=...=\gamma_{d}\}\:.
\]

\begin{prop}
\label{L2}Assume $m<d$ and that $\mathbf{G}$ is $m$-irreducible,
then there exists a unique $\mu_{F}\in\mathcal{M}(G_{d,m})$ with
$\mu_{F}=\sum_{\lambda\in\Lambda}p_{\lambda}\cdot A_{\lambda}^{-1}\mu_{F}$.
It also holds that $\dim_{H}\mu_{F}>0$.
\end{prop}

The measure $\mu_{F}$ is called the Furstenberg measure on $G_{d,m}$
corresponding to $\{A_{\lambda}^{-1}\}_{\lambda\in\Lambda}$ and $p$.
We can now state our main result:
\begin{thm}
\label{T3}If $m<d$, if $\mathbf{G}$ is $m$-irreducible, and if
\[
\dim_{H}^{*}\mu_{F}+D(\mu)>(m+1)(d-m),
\]
then $\pi_{v}\mu$ is exact dimensional with $\dim\pi_{v}\mu=D(\mu)$
for each $v\in\mathcal{V}$.
\end{thm}

\begin{rem}
As mentioned in the introduction, if $m=d$ then it follows from Theorem
2.6 in \cite{FH} that $\dim\pi_{v}\mu=D(\mu)$ for all $v\in\mathcal{V}$.
\end{rem}

\section{\label{S3}Proof of the main result}

For the remainder of this paper we assume $m<d$, $\mathbf{G}$ is
$m$-irreducible, and $\dim_{H}^{*}\mu_{F}+D(\mu)>(m+1)(d-m)$.

\subsection{\label{S1.2}Disintegration of measures}

For the proof of Theorem \ref{T3} we shall need to disintegrate the
measures $\mu$ and $\{\pi_{v}\mu\}_{v\in\mathcal{V}}$. We now define
these disintegrations and state some of their properties, for further
details see chapter 3 of \cite{FH}.

Let $\mathcal{B}$ be the Borel $\sigma$-algebra of $\mathbb{R}^{d}$,
let $X$ be a metric space, let $\theta\in\mathcal{M}(X)$, let $K$
be the support of $\theta$, and let $f:X\rightarrow\mathbb{R}^{d}$
be continuous. Then there exists a family $\{\theta_{x}\}_{x\in X}\subset\mathcal{M}(X)$,
which will be called the disintegration of $\theta$ with respect
to $f^{-1}\mathcal{B}$, such that:

\emph{(a)} For $\theta$-a.e. $x\in X$ the measure $\theta_{x}$
is supported on $K\cap f^{-1}(f(x))$.

\emph{(b) }For each $g\in L^{1}(\theta)$ and $\theta$-a.e. $x\in X$
we have
\[
\int g\:d\theta_{x}=\underset{\epsilon\downarrow0}{\lim}\frac{1}{f\theta(B(fx,\epsilon))}\cdot\int_{f^{-1}(B(fx,\epsilon))}g\:d\theta=\frac{d(f\theta^{g})}{d(f\theta)}(fx),
\]
where $\theta^{g}(E)=\int_{E}g\:d\theta$ for each Borel set $E\subset X$.
Here $\frac{d(f\theta^{g})}{d(f\theta)}$ stands for the Radon\textendash Nikodym
derivative of $f\theta^{g}$ with respect to $f\theta$.

\emph{(c) }For each $g\in L^{1}(\theta)$ the map that takes $x\in X$
to $\int g\:d\theta_{x}$ is $f^{-1}\mathcal{B}$ measurable and
\[
\int g\:d\theta_{x}=E_{\theta}[g\mid f^{-1}\mathcal{B}](x)\mbox{ for }\theta\mbox{-a.e. }x\in X\:.
\]
Here $E_{\theta}[g\mid f^{-1}\mathcal{B}]$ is the conditional expectation
of $g$ given $f^{-1}\mathcal{B}$ with respect to $\theta$.

We shall use the following notations for the disintegrations of $\mu$
and $\{\pi_{v}\mu\}_{v\in\mathcal{V}}$. For a subspace $W\subset\mathbb{R}^{d}$
set $\mathcal{B}_{W}=P_{W^{\perp}}^{-1}(\mathcal{B})$, and for $\theta\in\mathcal{M}(\mathbb{R}^{d})$
let $\{\theta_{W,x}\}_{x\in\mathbb{R}^{d}}$ be the disintegration
of $\theta$ with respect to $\mathcal{B}_{W}$. Given $v\in\mathbb{R}^{d|\Lambda|}$
set $\mathcal{F}_{v,W}=\pi_{v}^{-1}\circ P_{W^{\perp}}^{-1}(\mathcal{B})$
and let $\{\mu_{v,W,\omega}\}_{\omega\in\Omega}$ be the disintegration
of $\mu$ with respect to $\mathcal{F}_{v,W}$.

\subsection{Statement of auxiliary claims}

We now state some auxiliary claims which will be used in the proof
of Theorem \ref{T3}. The proofs are deferred to subsequent sections
in order to make the argument for Theorem \ref{T3} more transparent.
First we state Proposition \ref{P3} whose proof, which is given in
Section \ref{S5} below, requires ergodic theory and some results
from the random matrix theory presented in \cite{BL}.

Define $F:\mathcal{V}\rightarrow[0,\infty)$ by
\[
F(v)=-\frac{1}{\gamma_{d}}\cdot\int_{G_{d,m}}H_{\mu}(\mathcal{P}\mid\mathcal{F}_{v,W})\:d\mu_{F}(W)\;\mbox{ for }v\in\mathcal{V},
\]
where
\[
\mathcal{P}=\{\{\omega\in\Omega\::\:\omega_{0}=\lambda\}\in\mathcal{F}\::\:\lambda\in\Lambda\}
\]
and $H_{\mu}(\mathcal{P}\mid\mathcal{F}_{v,W})$ is the conditional
entropy of $\mathcal{P}$ given $\mathcal{F}_{v,W}$ with respect
to $\mu$.
\begin{prop}
\label{P3}For each $v\in\mathcal{V}$ and for $\mu\times\mu_{F}$-a.e.
$(\omega,W)\in\Omega\times G_{d,m}$ the measure $\pi_{v}\mu_{v,W,\omega}$
is exact dimensional with $\dim(\pi_{v}\mu_{v,W,\omega})=F(v)$.
\end{prop}

The rest of the auxiliary Lemmas will be proven is Section \ref{S6}.
\begin{lem}
\label{L6}Let $v\in\mathbb{R}^{d|\Lambda|}$ and $W\in G_{d,m}$,
then $(\pi_{v}\mu)_{W,\pi_{v}(\omega)}=\pi_{v}\mu_{v,W,\omega}$ for
$\mu$-a.e. $\omega\in\Omega$.
\end{lem}

The following semi-continuity lemma makes it possible to utilize Proposition
\ref{P3}.
\begin{lem}
\label{L7}The function $F$ is upper semi-continuous.
\end{lem}

\begin{lem}
\label{L8}For $v\in\mathcal{V}$ we have $\pi_{v}\mu\perp\mathcal{L}eb_{d}$,
where $\mathcal{L}eb_{d}$ is the Lebesgue measure of $\mathbb{R}^{d}$.
\end{lem}

The proof of the following lemma relies on results found in \cite{M2},
which are obtained by the use of Fourier analytic techniques. This
lemma makes it possible to use the assumption $\dim_{H}^{*}\mu_{F}+D(\mu)>(m+1)(d-m)$.
\begin{lem}
\label{L9}Let $\theta\in\mathcal{M}(\mathbb{R}^{d})$, let $1\le l<d$
be an integer, and set $s=\dim_{H}\theta$.

(a) If $s\le d-l$ then for $0\le t\le s$
\[
\dim_{H}\{W\in G_{d,l}\::\:\mathrm{essinf}_{\theta}\{\dim_{H}(\theta_{W,x})\::\:x\in\mathbb{R}^{d}\}>s-t\}\le(l-1)(d-l)+t\:.
\]

(b) If $s>d-l$ then for $s-l(d-l)\le t\le d-l$
\[
\dim_{H}\{W\in G_{d,l}\::\:\mathrm{essinf}_{\theta}\{\dim_{H}(\theta_{W,x})\::\:x\in\mathbb{R}^{d}\}>s-t\}\le l(d-l)+t-s\:.
\]

(c) If $s>d-l$ then
\[
\dim_{H}\{W\in G_{d,l}\::\:\mathrm{essinf}_{\theta}\{\dim_{H}(\theta_{W,x})\::\:x\in\mathbb{R}^{d}\}<s-d+l\}\le(l+1)(d-l)-s\:.
\]

\end{lem}

The proof for the following lemma is an adaptation of an argument
given in the proof of part (a) of Theorem 4.3 from \cite{JPS}.
\begin{lem}
\label{L10}For each $v\in\mathbb{R}^{d|\Lambda|}$ and for $\pi_{v}\mu$-a.e.
$x\in\mathbb{R}^{d}$ 
\[
\underset{\epsilon\downarrow0}{\limsup}\frac{\log\pi_{v}\mu(B(x,\epsilon))}{\log\epsilon}\le D(\mu)\:.
\]

\end{lem}

Let $\Lambda^{*}$ be the set of finite words over $\Lambda$. Given
a set of transformations (or matrices) $\{f_{\lambda}\}_{\lambda\in\Lambda}$,
that can be composed with one another, we set $f_{w}=f_{\lambda_{1}}\circ...\circ f_{\lambda_{k}}$
for $k\ge1$ and $\lambda_{1}\cdot...\cdot\lambda_{k}=w\in\Lambda^{*}$.
Given a set of real numbers $\{a_{\lambda}\}_{\lambda\in\Lambda}$
we set $a_{w}=a_{\lambda_{1}}\cdot...\cdot a_{\lambda_{k}}$. We also
set $f_{\emptyset}=Id$ and $a_{\emptyset}=1$, where $\emptyset\in\Lambda^{*}$
is the empty word.
\begin{lem}
\label{L11}Let $n\ge1$, let $\mathbf{G}'\subset Gl(d,\mathbb{R})$
be the closure of the group generated by $\{A_{w}\}_{w\in\Lambda^{n}}$,
set $p'=(p_{w})_{w\in\Lambda^{n}}$, set $\mu'=(p')^{\mathbb{N}}$,
and let $0>\gamma_{1}'\ge...\ge\gamma_{d}'>-\infty$ be the Lyapunov
exponents corresponding to $\mu'$ and $\{A_{w}\}_{w\in\Lambda^{n}}$.
Then $\mathbf{G}'$ is $m$-irreducible, $\gamma_{i}'=n\cdot\gamma_{i}$
for $1\le i\le d$, and $\mu_{F}'=\mu_{F}$, where $\mu_{F}'$ is
the Furstenberg measure corresponding to $\{A_{w}^{-1}\}_{w\in\Lambda^{n}}$
and $p'$ (see Proposition \ref{L2} above).
\end{lem}

\subsection{Proof of Theorem \ref{T3}}

By using Proposition \ref{P3} and Lemmas \ref{L6} to \ref{L11}
we shall now prove Theorem \ref{T3}.
\begin{lem}
\label{L12}If $\left\Vert A_{\lambda}\right\Vert <\frac{1}{2}$ for
each $\lambda\in\Lambda$, then $D(\mu)\in(d-m,d]$ and $F(v)\ge D(\mu)-d+m$
for each $v\in\mathcal{V}$.
\end{lem}

\textbf{\emph{Proof of Lemma \ref{L12}:}}\emph{ }Since $\mathcal{V}$
is non empty (by assumption) and since it is an open subset of $\mathbb{R}^{d|\Lambda|}$,
it follows that $\mathcal{L}eb_{d|\Lambda|}(\mathcal{V})>0$. From
part (b) of Theorem 1.9 in \cite{JPS} it follows that if $D(\mu)>d$,
then for $\mathcal{L}eb_{d|\Lambda|}$-a.e. $v\in\mathcal{V}$ we
have $\pi_{v}\mu\ll\mathcal{L}eb_{d}$. This together with Lemma \ref{L8}
shows that $D(\mu)\le d$. Since
\[
\dim_{H}^{*}\mu_{F}\le\dim_{H}G_{d,m}=m(d-m)
\]
 and
\[
\dim_{H}^{*}\mu_{F}+D(\mu)>m(d-m)+d-m,
\]
it follows that $D(\mu)\in(d-m,d]$. From this and from part (a) of
Theorem 1.9 in \cite{JPS} we get that $\dim_{H}\pi_{v}\mu=D(\mu)$
for $\mathcal{L}eb_{d|\Lambda|}$-a.e. $v\in\mathbb{R}^{d|\Lambda|}$.
Since $\mathcal{V}$ is open it follows that the set
\[
\mathcal{Q}=\{v\in\mathcal{V}\::\:\dim_{H}\pi_{v}\mu=D(\mu)\}
\]
is dense in $\mathcal{V}$.$\newline$Fix $v\in\mathcal{Q}$, then
from Proposition \ref{P3}, from Lemma \ref{L6}, and from (\ref{E2}),
it follows that for $\mu_{F}$-a.e. $W\in G_{d,m}$ we have for $\pi_{v}\mu$-a.e.
$x\in\mathbb{R}^{d}$ that $\dim_{H}(\pi_{v}\mu)_{W,x}=F(v)$. Set
\[
\mathcal{E}=\{W\in G_{d,m}\::\:\mathrm{essinf}_{\pi_{v}\mu}\{\dim_{H}(\pi_{v}\mu)_{W,x}\::\:x\in\mathbb{R}^{d}\}<D(\mu)-d+m\},
\]
then from $\dim_{H}\pi_{v}\mu=D(\mu)>d-m$ and from part (c) of Lemma
\ref{L9} we get 
\[
\dim_{H}(\mathcal{E})\le(m+1)(d-m)-D(\mu)\:.
\]
Since $\dim_{H}^{*}\mu_{F}>(m+1)(d-m)-D(\mu)$ it follows that $\mu_{F}(G_{d,m}\setminus\mathcal{E})>0$,
and so there exist $W\in G_{d,m}$ and $x\in\mathbb{R}^{d}$ with
\[
F(v)=\dim_{H}(\pi_{v}\mu)_{W,x}\ge D(\mu)-d+m\:.
\]
Since this holds for each $v\in\mathcal{Q}$ and since $\mathcal{Q}$
is dense in $\mathcal{V}$, it follows from Lemma \ref{L7}\emph{
}that $F(v)\ge D(\mu)-d+m$ for each $v\in\mathcal{V}$. $\square$

\textbf{\emph{Proof of Theorem \ref{T3}:}} Let $v\in\mathcal{V}$
be given. Assume first that $\left\Vert A_{\lambda}\right\Vert <\frac{1}{2}$
for each $\lambda\in\Lambda$, then from Lemma \ref{L12} we get $F(v)\ge D(\mu)-d+m\in(0,m]$.
From this, from Proposition \ref{P3}, and from Lemma \ref{L6} it
follows that
\begin{equation}
\dim_{H}(\pi_{v}\mu)_{W,x}\ge D(\mu)-d+m\mbox{ for }\pi_{v}\mu\times\mu_{F}\mbox{-a.e. }(x,W)\:.\label{E13}
\end{equation}
Set $s=\dim_{H}(\pi_{v}\mu)$. If $s<D(\mu)-d+m$ then clearly
\[
\mathrm{essinf}_{\pi_{v}\mu}\{\dim_{H}(\pi_{v}\mu)_{W,x}\::\:x\in\mathbb{R}^{d}\}<D(\mu)-d+m
\]
for each $W\in G_{d,m}$, and so we must have $s\ge D(\mu)-d+m$.
Assume by contradiction that $D(\mu)-d+m\le s<D(\mu)$, let
\[
0<\epsilon<\min\left\{ \begin{array}{c}
D(\mu)-d+m,\\
D(\mu)-s,\\
\dim_{H}^{*}\mu_{F}+D(\mu)-(m+1)(d-m)
\end{array}\right\} ,
\]
set
\[
t=\begin{cases}
\min\{2(d-m)-D(\mu)+\epsilon,s\} & ,\mbox{ if }s\le d-m\\
d-m+s-D(\mu)+\epsilon & ,\mbox{ if }s>d-m
\end{cases},
\]
and set
\[
\mathcal{E}=\{W\in G_{d,m}\::\:\mathrm{essinf}_{\pi_{v}\mu}\{\dim_{H}(\pi_{v}\mu)_{W,x}\::\:x\in\mathbb{R}^{d}\}>s-t\}\:.
\]
If $s\le d-m$ then
\[
D(\mu)-d+m\le s\le d-m,
\]
so $0\le t\le s$, and so from part (a) of Lemma \ref{L9}
\[
\dim_{H}(\mathcal{E})\le(m-1)(d-m)+t\le(m+1)(d-m)-D(\mu)+\epsilon<\dim_{H}^{*}\mu_{F}\:.
\]
If $s>d-m$ then
\[
t-(s-m(d-m))>d-m-D(\mu)+m(d-m)\ge m(d-m)-m\ge0
\]
and
\[
d-m-t=D(\mu)-s-\epsilon>0,
\]
so $s-m(d-m)\le t\le d-m$, and so from part (b) of Lemma \ref{L9}
\[
\dim_{H}(\mathcal{E})\le m(d-m)+t-s=(m+1)(d-m)-D(\mu)+\epsilon<\dim_{H}^{*}\mu_{F}\:.
\]
In any case we have $\dim_{H}(\mathcal{E})<\dim_{H}^{*}\mu_{F}$,
so $\mu_{F}(G_{d,m}\setminus\mathcal{E})>0$, and so
\[
\pi_{v}\mu\times\mu_{F}\{(x,W)\::\:\dim(\pi_{v}\mu)_{W,x}\le s-t+\frac{\epsilon}{2}\}>0\:.
\]
But this gives a contradiction to (\ref{E13}) since if $s\le d-m$
then
\begin{multline*}
s-t+\frac{\epsilon}{2}=\max\{s-(2(d-m)-D(\mu)+\epsilon),0\}+\frac{\epsilon}{2}\\
\le\max\{D(\mu)-d+m-\epsilon,0\}+\frac{\epsilon}{2}=D(\mu)-d+m-\frac{\epsilon}{2},
\end{multline*}
and if $s>d-m$ then 
\[
s-t+\frac{\epsilon}{2}=D(\mu)-d+m-\frac{\epsilon}{2}\:.
\]
It follows that we must have $\dim_{H}(\pi_{v}\mu)=s\ge D(\mu)$,
and so from Lemma \ref{L10} and (\ref{E2}) we obtain that $\pi_{v}\mu$
is exact dimensional with $\dim\pi_{v}\mu=D(\mu)$. This proves the
theorem if $\left\Vert A_{\lambda}\right\Vert <\frac{1}{2}$ for each
$\lambda\in\Lambda$.$\newline$Now we prove the general case. Let
$n\ge1$ be such that $\left\Vert A_{w}\right\Vert <\frac{1}{2}$
for each $w\in\Lambda^{n}$. Since the SSC holds for $\{\varphi_{v,\lambda}\}_{\lambda\in\Lambda}$
it clearly holds for $\{\varphi_{v,w}\}_{w\in\Lambda^{n}}$. For $\omega\in(\Lambda^{n})^{\mathbb{N}}$
set $\pi_{v}'(\omega)=\underset{n}{\lim}\:\varphi_{v,\omega_{0}}\circ...\circ\varphi_{v,\omega_{n}}(0)$,
set $p'=(p_{w})_{w\in\Lambda^{n}}$, set $\mu'=(p')^{\mathbb{N}}$,
let $0>\gamma_{1}'\ge...\ge\gamma_{d}'>-\infty$ be the Lyapunov exponents
corresponding to $\mu'$ and $\{A_{w}\}_{w\in\Lambda^{n}}$, and let
$\mathbf{G}'\subset Gl(d,\mathbb{R})$ be the closure of the group
generated by $\{A_{w}\}_{w\in\Lambda^{n}}$. From Lemma \ref{L11}
we get that $\mathbf{G}'$ is $m$-irreducible, $\gamma_{i}'=n\cdot\gamma_{i}$
for $1\le i\le d$, and $\mu_{F}'=\mu_{F}$, where $\mu_{F}'$ is
the Furstenberg measure corresponding to $\{A_{w}^{-1}\}_{w\in\Lambda^{n}}$
and $p'$. Let $h_{\mu'}$ be the entropy of $\mu'$ (i.e. $h_{\mu'}=\sum_{w\in\Lambda^{n}}-p_{w}\cdot\log p_{w}$),
and let $D(\mu')$ be the Lyapunov dimension of $\mu'$ with respect
to the family $\{A_{w}\}_{w\in\Lambda^{n}}$ (see the definition in
Section \ref{S2} above). Since $h_{\mu'}=n\cdot h_{\mu}$ it follows
from the definition of the Lyapunov dimension that $D(\mu')=D(\mu)$,
hence
\[
\dim_{H}^{*}\mu_{F}'+D(\mu')=\dim_{H}^{*}\mu_{F}+D(\mu)>(m+1)(d-m)\:.
\]
Now from the first part of the proof we get that $\pi_{v}'\mu'$ is
exact dimensional with $\dim\pi_{v}'\mu'=D(\mu')=D(\mu)$. This completes
the proof since $\pi_{v}\mu=\pi_{v}'\mu'$. $\square$

\section{\label{S4}Auxiliary results from the theory of random matrices}

In this section we translate results found in \cite{BL} to suit our
needs. These results will be used in the proofs of Propositions \ref{L2}
and \ref{P3}.
\begin{defn}
\label{D14}Given $q\ge2$, $1\le l<q$, and $\mathbf{S}\subset Gl(q,\mathbb{R})$,
we say that $\mathbf{S}$ is $l$-strongly irreducible if there does
not exist a finite family of proper linear subspaces $W_{1},...,W_{k}$
of $\mathcal{A}^{l}(\mathbb{R}^{q})$ with
\[
\mathcal{A}^{l}M(W_{1}\cup...\cup W_{k})=W_{1}\cup...\cup W_{k}\mbox{ for each }M\in\mathbf{S}\:.
\]
When $l=1$ we say that $\mathbf{S}$ is strongly irreducible.
\end{defn}

\begin{rem}
\label{R14}Given $q\ge2$, $1\le l<q$, and linear subspaces $W_{1},...,W_{k}$
of $\mathcal{A}^{l}(\mathbb{R}^{q})$, the set
\[
\{M\in Gl(q,\mathbb{R})\::\:\mathcal{A}^{l}M(W_{1}\cup...\cup W_{k})=W_{1}\cup...\cup W_{k}\}
\]
is a closed subgroup of $Gl(q,\mathbb{R})$.
\end{rem}

\begin{defn}
Given $q\ge2$, $1\le l<q$, and $\mathbf{S}\subset Gl(q,\mathbb{R})$,
we say that $\mathbf{S}$ is $l$-contracting if there exists a sequence
$\{M_{n}\}_{n=1}^{\infty}\subset\mathbf{S}$ such that
\[
\{\left\Vert \mathcal{A}^{l}M_{n}\right\Vert ^{-1}\cdot\mathcal{A}^{l}M_{n}\::\:n\ge1\}
\]
converges to a rank-one matrix. When $l=1$ we say that $\mathbf{S}$
is contracting.
\end{defn}

Throughout this section $\mathbf{T}\subset Gl(d,\mathbb{R})$ will
denote the closure of the semigroup generated by $\{A_{\lambda}^{-1}\}_{\lambda\in\Lambda}$.
Let $q\ge1$ be the dimension of $\mathcal{A}^{m}(\mathbb{R}^{d})$,
then given $M\in Gl(d,\mathbb{R})$ we may view $\mathcal{A}^{m}M$
as a member of $Gl(q,\mathbb{R})$. Let $\widetilde{\mathbf{T}}\subset Gl(q,\mathbb{R})$
be the closure of the semigroup generated by $\{\mathcal{A}^{m}A_{\lambda}^{-1}\}_{\lambda\in\Lambda}$.
Recall that we assume $m<d$ and $\mathbf{G}$ is $m$-irreducible.
\begin{lem}
\label{L15}$\widetilde{\mathbf{T}}$ is contracting and strongly
irreducible, and $\mathbf{T}$ is $m$-contracting and $m$-strongly
irreducible.
\end{lem}

\textbf{\emph{Proof of Lemma \ref{L15}:}} Since $\mathbf{G}$ is
$m$-irreducible it follows from remark \ref{R14} that $\{A_{\lambda}^{-1}\}_{\lambda\in\Lambda}$
is $m$-irreducible, and so $\widetilde{\mathbf{T}}$ is irreducible.
Let $\infty>\gamma_{1}'\ge...\ge\gamma_{d}'>0$ be the Lyapunov exponents
corresponding to $\mu$ and $\{A_{\lambda}^{-1}\}_{\lambda\in\Lambda}$,
then $\gamma_{i}'=-\gamma_{d-i+1}$ for $1\le i\le d$. Let $\eta_{1}\ge\eta_{2}$
be the the two upper Lyapunov exponents corresponding to $\mu$ and
$\{\mathcal{A}^{m}A_{\lambda}^{-1}\}_{\lambda\in\Lambda}$. From an
argument given in the proof of Theorem IV.1.2 in \cite{BL} we get
\[
\eta_{1}=\sum_{i=1}^{m}\gamma_{i}'\;\mbox{ and }\;\eta_{2}=\sum_{i=1}^{m-1}\gamma_{i}'+\gamma_{m+1}',
\]
hence from the definition of $m$
\[
\eta_{1}=\sum_{i=1}^{m}\gamma_{i}'=-\sum_{i=1}^{m}\gamma_{d-i+1}>-\sum_{i=1}^{m-1}\gamma_{d-i+1}-\gamma_{d-m}=\sum_{i=1}^{m-1}\gamma_{i}'+\gamma_{m+1}'=\eta_{2}\:.
\]
From this, from the irreducibility of $\widetilde{\mathbf{T}}$, and
from Theorem III.6.1 in \cite{BL}, we get that $\widetilde{\mathbf{T}}$
is contracting and strongly irreducible. From this and remark \ref{R14}
it follows that $\{\mathcal{A}^{m}A_{\lambda}^{-1}\}_{\lambda\in\Lambda}$
is strongly irreducible, and so $\mathbf{T}$ is $m$-strongly irreducible.
Since $\widetilde{\mathbf{T}}$ is contracting and since $\{\mathcal{A}^{m}A_{w}^{-1}\::\:w\in\Lambda^{*}\}$
is dense in $\widetilde{\mathbf{T}}$, it follows that $\{\mathcal{A}^{m}A_{w}^{-1}\::\:w\in\Lambda^{*}\}$
is contracting. This shows that $\mathbf{T}$ is $m$-contracting.
$\square$

Let $\left\langle \cdot,\cdot\right\rangle $ be the usual scalar
product on $\mathbb{R}^{d}$. As in Section III.5 of \cite{BL} we
define a scalar product on $\mathcal{A}^{m}(\mathbb{R}^{d})$ by the
formula
\[
\left\langle x_{1}\wedge...\wedge x_{m},y_{1}\wedge...\wedge y_{m}\right\rangle =\det\left[\{\left\langle x_{i},y_{j}\right\rangle \}_{i,j=1}^{m}\right]\:.
\]
Let $P(\mathcal{A}^{m}(\mathbb{R}^{d}))$ be the projective space
of $\mathcal{A}^{m}(\mathbb{R}^{d})$. Given $\bar{\xi},\bar{\eta}\in P(\mathcal{A}^{m}(\mathbb{R}^{d}))$
set 
\[
d_{P(\mathcal{A}^{m}(\mathbb{R}^{d}))}(\bar{\xi},\bar{\eta})=\left(1-\left\langle \xi,\eta\right\rangle ^{2}\right)^{1/2},
\]
where $\xi$ and $\eta$ are unit vectors in $\mathcal{A}^{m}(\mathbb{R}^{d})$
with directions $\bar{\xi}$ and $\bar{\eta}$. As shown in Section
III.4 of \cite{BL}, $d_{P(\mathcal{A}^{m}(\mathbb{R}^{d}))}$ is
a metric on $P(\mathcal{A}^{m}(\mathbb{R}^{d}))$.$\newline$Given
independent sets $\{x_{1},...,x_{m}\},\{y_{1},...,y_{m}\}\subset\mathbb{R}^{d}$,
there exists a constant $a\in\mathbb{R}\setminus\{0\}$ with
\[
y_{1}\wedge...\wedge y_{m}=a\cdot x_{1}\wedge...\wedge x_{m}
\]
if and only if
\[
\mathrm{span}\{y_{1},...,y_{m}\}=\mathrm{span}\{x_{1},...,x_{m}\}\:.
\]
Define a map $\psi:G_{d,m}\rightarrow P(\mathcal{A}^{m}(\mathbb{R}^{d}))$
by 
\[
\psi(W)=\mathbb{R}\cdot x_{1}\wedge...\wedge x_{m}\;\mbox{ if }\mathrm{span}\{x_{1},...,x_{m}\}=W\in G_{d,m}\:.
\]
It is not hard to check that there exists a constant $C\in(1,\infty)$
with
\begin{equation}
C^{-1}\cdot d_{G_{d,m}}(W,U)\le\left(d_{P(\mathcal{A}^{m}(\mathbb{R}^{d}))}(\psi(W),\psi(U))\right)^{2}\le C\cdot d_{G_{d,m}}(W,U)\label{E21}
\end{equation}
for all $W,U\in G_{d,m}$, where $d_{G_{d,m}}$ is the metric defined
above in Section \ref{S2}. Hence $\psi$ is an embedding of $G_{d,m}$
into $P(\mathcal{A}^{m}(\mathbb{R}^{d}))$. Now we can prove Proposition
\ref{L2}\emph{.}

\textbf{\emph{Proof of Proposition \ref{L2}:}} From Lemma \ref{L15}
and Theorem IV.1.2 in \cite{BL} it follows that there exists a unique
$\theta\in\mathcal{M}(P(\mathcal{A}^{m}(\mathbb{R}^{d})))$ with $\theta=\sum_{\lambda\in\Lambda}p_{\lambda}\cdot\mathcal{A}^{m}A_{\lambda}^{-1}\theta$.
Since $\psi(G_{d,m})$ is compact and $\mathcal{A}^{m}M(\psi(G_{d,m}))=\psi(G_{d,m})$
for each $M\in Gl(d,\mathbb{R})$, it follows from Lemma I.3.5 in
\cite{BL} that there exits $\theta'\in\mathcal{M}(\psi(G_{d,m}))$
with $\theta'=\sum_{\lambda\in\Lambda}p_{\lambda}\cdot\mathcal{A}^{m}A_{\lambda}^{-1}\theta'$.
By the uniqueness of $\theta$ it follows that $\theta=\theta'$,
and so $\theta$ is supported on $\psi(G_{d,m})$.$\newline$Set $\mu_{F}=\psi^{-1}\theta$,
then
\[
\mu_{F}=\psi^{-1}\theta=\sum_{\lambda\in\Lambda}p_{\lambda}\cdot\psi^{-1}\circ\mathcal{A}^{m}A_{\lambda}^{-1}\theta=\sum_{\lambda\in\Lambda}p_{\lambda}\cdot A_{\lambda}^{-1}\circ\psi^{-1}\theta=\sum_{\lambda\in\Lambda}p_{\lambda}\cdot A_{\lambda}^{-1}\mu_{F}\:.
\]
Since $\psi$ is an embedding the uniqueness of $\mu_{F}$ follows
from the uniqueness of $\theta$. From Corollary VI.4.2 in \cite{BL}
and the remarks following it it follows that $\dim_{H}\theta>0$.
From this and from (\ref{E21}) we obtain $\dim_{H}\mu_{F}>0$. This
completes the proof of the Lemma. $\square$

Given $a_{1},...,a_{d}\in\mathbb{R}$ let $diag(a_{1},...,a_{d})$
denote the $d\times d$ matrix $D$ with
\[
D_{i,j}=\begin{cases}
a_{i} & ,\mbox{ if }i=j\\
0 & ,\mbox{ if }i\ne j
\end{cases}\mbox{ for }1\le i,j\le d\:.
\]
Given $M\in Gl(d,\mathbb{R})$ there exist orthogonal matrices $U,V\in Gl(d,\mathbb{R})$
with $M=UDV$, where $D=diag(\alpha_{1}(M),...,\alpha_{d}(M))$. We
call the product $UDV$ a singular value decomposition of $M$. Note
that $V^{*}e_{i}$ is an eigenvector of $M^{*}M$ with eigenvalue
$\alpha_{i}(M)^{2}$ for each $1\le i\le d$. Here $\{e_{i}\}_{i=1}^{d}$
is the standard basis of $\mathbb{R}^{d}$ and $M^{*}$ is the transpose
of $M$.
\begin{lem}
\label{L16}For each $\omega\in\Omega$ and $n\ge1$ set $D_{n,\omega}=diag(\alpha_{1}(A_{\omega|_{n}}),...,\alpha_{d}(A_{\omega|_{n}}))$,
let $U_{n,\omega}D_{n,\omega}V_{n,\omega}$ be a singular value decomposition
of $A_{\omega|_{n}}$, and set $W_{n}(\omega)=span\{U_{n,\omega}e_{d-m+1},...,U_{n,\omega}e_{d}\}$.
Then for $\mu$-a.e. $\omega\in\Omega$ there exists $W(\omega)\in G_{d,m}$
such that $\{W_{n}(\omega)\}_{n=1}^{\infty}$ converges to $W(\omega)$
in $G_{d,m}$.
\end{lem}

\textbf{\emph{Proof of Lemma \ref{L16}:}} From Lemma \ref{L15} we
get that $\widetilde{\mathbf{T}}$ is a contracting and strongly irreducible
subset of $Gl(q,\mathbb{R})$. Hence we may apply proposition III.3.2
in \cite{BL} on the i.i.d. sequence $\{\mathcal{A}^{m}A_{\omega_{n}}^{-1}\}_{n=0}^{\infty}$.
For each $\omega\in\Omega$ and $n\ge1$ set $M_{n,\omega}=A_{\omega_{n-1}}^{-1}\cdot...\cdot A_{\omega_{0}}^{-1}$,
set $\xi_{n,\omega}=U_{n,\omega}e_{d-m+1}\wedge...\wedge U_{n,\omega}e_{d}$,
and set
\[
\widetilde{W}_{n}(\omega)=\{\eta\in\mathcal{A}^{m}(\mathbb{R}^{d})\::\:\mathcal{A}^{m}M_{n.\omega}^{*}M_{n.\omega}\eta=\alpha_{1}(\mathcal{A}^{m}M_{n.\omega}^{*}M_{n.\omega})\cdot\eta\}\:.
\]
From part (b) of proposition III.3.2 it follows that for $\mu$-a.e.
$\omega\in\Omega$
\[
\alpha_{1}(\mathcal{A}^{m}M_{n.\omega}^{*}M_{n.\omega})>\alpha_{2}(\mathcal{A}^{m}M_{n.\omega}^{*}M_{n.\omega})
\]
for all $n$ large enough, and so $\widetilde{W}_{n}(\omega)$ is
$1$-dimensional for all $n$ large enough. From part (a) of proposition
III.3.2 it follows that for $\mu$-a.e. $\omega\in\Omega$ the sequence
$\{\widetilde{W}_{n}(\omega)\}_{n=1}^{\infty}$ converges to some
element in $P(\mathcal{A}^{m}(\mathbb{R}^{d}))$. For each $\omega\in\Omega$
and $n\ge1$ we have 
\begin{multline*}
M_{n.\omega}^{*}M_{n.\omega}U_{n,\omega}=(A_{\omega|_{n}}^{-1})^{*}A_{\omega|_{n}}^{-1}U_{n,\omega}\\
=(V_{n,\omega}^{-1}D_{n,\omega}^{-1}U_{n,\omega}^{-1})^{*}(V_{n,\omega}^{-1}D_{n,\omega}^{-1}U_{n,\omega}^{-1})U_{n,\omega}=U_{n,\omega}D_{n,\omega}^{-2},
\end{multline*}
and also from Lemma 5.3 in \cite{BL}
\begin{multline*}
\alpha_{1}(\mathcal{A}^{m}M_{n.\omega}^{*}M_{n.\omega})=\prod_{i=1}^{m}\alpha_{i}(M_{n.\omega}^{*}M_{n.\omega})=\prod_{i=1}^{m}\alpha_{i}(M_{n.\omega})^{2}\\
=\prod_{i=1}^{m}\alpha_{i}(A_{\omega|_{n}}^{-1})^{2}=\prod_{i=1}^{m}\alpha_{d-i+1}(A_{\omega|_{n}})^{-2}\:.
\end{multline*}
It follows that
\begin{multline*}
\mathcal{A}^{m}M_{n.\omega}^{*}M_{n.\omega}(\xi_{n,\omega})=U_{n,\omega}D_{n,\omega}^{-2}e_{d-m+1}\wedge...\wedge U_{n,\omega}D_{n,\omega}^{-2}e_{d}\\
=\prod_{i=1}^{m}\alpha_{d-i+1}(A_{\omega|_{n}})^{-2}\cdot\xi_{n,\omega}=\alpha_{1}(\mathcal{A}^{m}M_{n.\omega}^{*}M_{n.\omega})\cdot\xi_{n,\omega},
\end{multline*}
hence $\xi_{n,\omega}\in\widetilde{W}_{n}(\omega)$, and so for $\mu$-a.e.
$\omega\in\Omega$ we have $\mathbb{R}\cdot\xi_{n,\omega}=\widetilde{W}_{n}(\omega)$
for all $n$ large enough. This shows that for $\mu$-a.e. $\omega\in\Omega$
the sequence $\{\mathbb{R}\cdot\xi_{n,\omega}\}_{n=1}^{\infty}$ converges
in $P(\mathcal{A}^{m}(\mathbb{R}^{d}))$. Now since $\{\mathbb{R}\cdot\xi_{n,\omega}\}_{n=1}^{\infty}\subset\psi(G_{d,m})$,
since $\psi(G_{d,m})$ is compact, and since $\psi$ is an embedding,
it follows that
\[
\{W_{n}(\omega)\}_{n=1}^{\infty}=\{\psi^{-1}(\mathbb{R}\cdot\xi_{n,\omega})\}_{n=1}^{\infty}
\]
converges to some $W(\omega)$ in $G_{d,m}$. This completes the proof
of the lemma. $\square$
\begin{lem}
\label{L17}Let $U\in G_{d,m}$ be given and set
\[
\mathcal{S}_{U}=\{W\in G_{d,m}\::\:U^{\perp}+W\ne\mathbb{R}^{d}\},
\]
then $\mu_{F}(\mathcal{S}_{U})=0$.
\end{lem}

\textbf{\emph{Proof of Lemma \ref{L17}:}} Set $\theta=\psi\mu_{F}$,
then $\theta\in\mathcal{M}(P(\mathcal{A}^{m}(\mathbb{R}^{d})))$ and
\[
\theta=\sum_{\lambda\in\Lambda}p_{\lambda}\cdot\mathcal{A}^{m}A_{\lambda}^{-1}\theta\:.
\]
From the strong irreducibility of $\widetilde{\mathbf{T}}$ and from
proposition III.2.3 in \cite{BL}, it follows that
\[
\theta\{\mathbb{R}\cdot z\::\:z\in\mathcal{Q}\setminus\{0\}\}=0
\]
for every proper subspace $\mathcal{Q}$ of $\mathcal{A}^{m}(\mathbb{R}^{d})$.
Let $\{x_{1},...,x_{d-m}\}$ be a basis for $U^{\perp}$, set $\xi=x_{1}\wedge...\wedge x_{d-m}$,
and set 
\[
\mathcal{Q}=\{z\in\mathcal{A}^{m}(\mathbb{R}^{d})\::\:\xi\wedge z=0\},
\]
then $\mathcal{Q}$ is a proper subspace of $\mathcal{A}^{m}(\mathbb{R}^{d})$.
Now since
\begin{multline*}
\mu_{F}(\mathcal{S}_{U})=\mu_{F}\{W\in G_{d,m}\::\:\xi\wedge w_{1}\wedge...\wedge w_{m}=0\mbox{ where }\{w_{1},...,w_{m}\}\mbox{ is a basis for }W\}\\
=\mu_{F}\{W\in G_{d,m}\::\:\psi(W)=\mathbb{R}\cdot z\mbox{ where }z\in\mathcal{A}^{m}(\mathbb{R}^{d})\mbox{ and }\xi\wedge z=0\}\\
=\theta\{\mathbb{R}\cdot z\::\:z\in\mathcal{Q}\setminus\{0\}\}=0
\end{multline*}
the lemma follows. $\square$

\section{\label{S5}Proof of Proposition \ref{P3}}

Fix some $v\in\mathcal{V}$ and set $\pi=\pi_{v}$, $K=K_{v}$, $\varphi_{\lambda}=\varphi_{v,\lambda}$
for $\lambda\in\Lambda$, and $\mathcal{F}_{W}=\mathcal{F}_{v,W}$
and $\{\mu_{W,\omega}\}_{\omega\in\Omega}=\{\mu_{v,W,\omega}\}_{\omega\in\Omega}$
for $W\in G_{d,m}$. For $k\ge1$ and $\lambda_{0}\cdot...\cdot\lambda_{k-1}=w\in\Lambda^{k}$
let
\[
[w]=\{\omega\in\Omega\::\:\omega_{i}=\lambda_{i}\mbox{ for }0\le i<k\},
\]
and let $[\emptyset]=\Omega$. Given $\omega\in\Omega$ and $k\ge1$
set $\omega|_{k}=\omega_{0}\cdot...\cdot\omega_{k-1}\in\Lambda^{k}$
and $\omega|_{0}=\emptyset$.

In the proof of Proposition \ref{P3} we shall make use of the following
dynamical system. Let $\sigma:\Omega\rightarrow\Omega$ be the left
shift, i.e. $(\sigma\omega)_{k}=\omega_{k+1}$ for $\omega\in\Omega$
and $k\ge0$. Set $X=\Omega\times G_{d,m}$, for each $(\omega,W)\in X$
set $T(\omega,W)=(\sigma(\omega),A_{\omega_{0}}^{-1}\cdot W)$, and
set $\nu=\mu\times\mu_{F}$. Since $\mu_{F}$ is the unique member
in $\mathcal{M}(G_{d,m})$ with $\mu_{F}=\sum_{\lambda\in\Lambda}p_{\lambda}\cdot A_{\lambda}^{-1}\mu_{F}$,
it follows from Proposition 1.14 in \cite{BQ} that $(X,T,\nu)$ is
measure preserving and ergodic.
\begin{lem}
\label{L18}Let $E\subset\Omega$ be a Borel set, let $M\in Gl(d,\mathbb{R})$,
let $W\in G_{d,m}$, and set $\widetilde{B}=P_{W^{\perp}}\circ M(B(0,1))$,
then for $\mu$-a.e. $\omega\in\Omega$ 
\[
\mu_{W,\omega}(E)=\underset{\delta\downarrow0}{\lim}\frac{\mu(\pi^{-1}\circ P_{W^{\perp}}^{-1}(P_{W^{\perp}}\circ\pi(\omega)+\delta\cdot\widetilde{B})\cap E)}{\mu(\pi^{-1}\circ P_{W^{\perp}}^{-1}(P_{W^{\perp}}\circ\pi(\omega)+\delta\cdot\widetilde{B}))}\:.
\]

\end{lem}

\textbf{\emph{Proof of Lemma \ref{L18}:}} Let $\mu|_{E}$ be the
restriction of $\mu$ to $E$, i.e. $\mu|_{E}(F)=\mu(F\cap E)$ for
$F\in\mathcal{F}$. For $x\in W^{\perp}$ set $\left\Vert x\right\Vert _{\widetilde{B}}=\inf\{t>0\::\:t^{-1}\cdot x\in\widetilde{B}\}$,
i.e. $\left\Vert \cdot\right\Vert _{\widetilde{B}}$ is the Minkowski
functional corresponding to the convex and balanced set $\widetilde{B}$.
Clearly $\left\Vert \cdot\right\Vert _{\widetilde{B}}$ is a norm
on $W^{\perp}$, and 
\[
\delta\cdot\widetilde{B}=\{x\in W^{\perp}\::\:\left\Vert x\right\Vert _{\widetilde{B}}\le\delta\}\:\mbox{ for }\delta>0\:.
\]
Now from Theorem 4.2 in \cite{BL1} and the discussion preceding it,
and from property (b) in Section \ref{S1.2} above, we get that for
$\mu$-a.e. $\omega\in\Omega$
\begin{multline*}
\underset{\delta\downarrow0}{\lim}\frac{\mu(\pi^{-1}\circ P_{W^{\perp}}^{-1}(P_{W^{\perp}}\circ\pi(\omega)+\delta\cdot\widetilde{B})\cap E)}{\mu(\pi^{-1}\circ P_{W^{\perp}}^{-1}(P_{W^{\perp}}\circ\pi(\omega)+\delta\cdot\widetilde{B}))}\\
=\underset{\delta\downarrow0}{\lim}\frac{P_{W^{\perp}}\pi\mu|_{E}(P_{W^{\perp}}\circ\pi(\omega)+\delta\cdot\widetilde{B})}{P_{W^{\perp}}\pi\mu(P_{W^{\perp}}\circ\pi(\omega)+\delta\cdot\widetilde{B})}\\
=\frac{dP_{W^{\perp}}\pi\mu|_{E}}{dP_{W^{\perp}}\pi\mu}(P_{W^{\perp}}\circ\pi(\omega))=\mu_{W,\omega}(E),
\end{multline*}
which proves the Lemma. $\square$
\begin{lem}
\label{L19}For each $W\in G_{d,m}$ and $k\ge0$ 
\[
\frac{\mu_{W,\omega}[\omega|_{k+1}]}{\mu_{W,\omega}[\omega|_{k}]}=\mu_{(A_{\omega|_{k}})^{-1}\cdot W,\sigma^{k}\omega}[\omega_{k}]\;\mbox{ for }\mu\mbox{-a.e. }\omega\in\Omega\:.
\]

\end{lem}

\textbf{\emph{Proof of Lemma \ref{L19}:}} For each $\lambda\in\Lambda$
and $\omega\in\Omega$ set $f_{\lambda}(\omega)=\lambda\cdot\omega$,
i.e. $f_{\lambda}(\omega)$ is the concatenation of $\lambda$ with
$\omega$. Let $W\in G_{d,m}$, $k\ge0$, and $w\in\Lambda^{k}$ be
given, and set $U=(A_{w})^{-1}\cdot W$. From property (b) stated
in Section \ref{S1.2} above and since $\mu(f_{w}(E))=p_{w}\cdot\mu(E)$
for each $E\in\mathcal{F}$, it follows that for $\mu$-a.e. $\omega\in\Omega$
\begin{multline}
\mu_{U,\sigma^{k}\omega}[\omega_{k}]=\underset{\delta\downarrow0}{\lim}\frac{\mu(\pi^{-1}\circ P_{U^{\perp}}^{-1}(B(P_{U^{\perp}}\circ\pi\circ\sigma^{k}(\omega),\delta))\cap[(\sigma^{k}\omega)|_{1}])}{\mu(\pi^{-1}\circ P_{U^{\perp}}^{-1}(B(P_{U^{\perp}}\circ\pi\circ\sigma^{k}(\omega),\delta)))}\\
\underset{\delta\downarrow0}{\lim}\frac{\mu(f_{w}(\pi^{-1}\circ P_{U^{\perp}}^{-1}(B(P_{U^{\perp}}\circ\pi\circ\sigma^{k}(\omega),\delta))\cap[(\sigma^{k}\omega)|_{1}]))}{\mu(f_{w}\circ\pi^{-1}\circ P_{U^{\perp}}^{-1}(B(P_{U^{\perp}}\circ\pi\circ\sigma^{k}(\omega),\delta)))}\:.\label{E12}
\end{multline}
Fix $\omega\in[w]$ and $\delta>0$, and set $\widetilde{B}=P_{W^{\perp}}\circ A_{w}(B(0,1))$.
Since $f_{w}\circ\pi^{-1}(x)=\pi^{-1}\circ\varphi_{w}(x)$ for $x\in K$,
\begin{multline}
f_{w}\circ\pi^{-1}\circ P_{U^{\perp}}^{-1}(B(P_{U^{\perp}}\circ\pi\circ\sigma^{k}(\omega),\delta))\\
=\pi^{-1}\circ\varphi_{w}(K\cap P_{U^{\perp}}^{-1}(B(P_{U^{\perp}}\circ\pi\circ\sigma^{k}(\omega),\delta)))\\
=\pi^{-1}\circ\varphi_{w}(K)\cap\pi^{-1}\circ\varphi_{w}(\pi\circ\sigma^{k}(\omega)+U+B(0,\delta))\\
=[\omega|_{k}]\cap\pi^{-1}\circ\varphi_{w}(\pi\circ\sigma^{k}(\omega)+U+B(0,\delta))\:.\label{E11}
\end{multline}
From $\varphi_{w}\circ\pi=\pi\circ f_{w}$ and $\omega|_{k}=w$ we
get
\begin{multline*}
\varphi_{w}(\pi\circ\sigma^{k}(\omega)+U+B(0,\delta))\\
=\pi\circ f_{w}\circ\sigma^{k}(\omega)+A_{w}\cdot U+A_{w}(B(0,\delta))\\
=\pi(\omega)+W+\delta\cdot A_{w}(B(0,1))\\
=W+P_{W^{\perp}}\circ\pi(\omega)+\delta\cdot P_{W^{\perp}}(A_{w}(B(0,1)))\\
=P_{W^{\perp}}^{-1}(P_{W^{\perp}}\circ\pi(\omega)+\delta\cdot\widetilde{B})\:.
\end{multline*}
From this and from (\ref{E11}) we obtain
\[
f_{w}\circ\pi^{-1}\circ P_{U^{\perp}}^{-1}(B(P_{U^{\perp}}\circ\pi\circ\sigma^{k}(\omega),\delta))=[\omega|_{k}]\cap\pi^{-1}\circ P_{W^{\perp}}^{-1}(P_{W^{\perp}}\circ\pi(\omega)+\delta\cdot\widetilde{B}),
\]
for each $\omega\in[w]$ and $\delta>0$. It now follows from (\ref{E12})
and Lemma \ref{L18} that for $\mu$-a.e. $\omega\in[w]$
\begin{multline*}
\mu_{U,\sigma^{k}\omega}[\omega_{k}]=\underset{\delta\downarrow0}{\lim}\frac{\mu([\omega|_{k}]\cap\pi^{-1}\circ P_{W^{\perp}}^{-1}(P_{W^{\perp}}\circ\pi(\omega)+\delta\cdot\widetilde{B})\cap f_{w}([(\sigma^{k}\omega)|_{1}]))}{\mu([\omega|_{k}]\cap\pi^{-1}\circ P_{W^{\perp}}^{-1}(P_{W^{\perp}}\circ\pi(\omega)+\delta\cdot\widetilde{B}))}\\
=\frac{\mu_{W,\omega}([\omega|_{k}]\cap f_{w}[(\sigma^{k}\omega)|_{1}])}{\mu_{W,\omega}[\omega|_{k}]}=\frac{\mu_{W,\omega}[\omega|_{k+1}]}{\mu_{W,\omega}[\omega|_{k}]}\:.
\end{multline*}
This proves the lemma since $U=(A_{\omega|_{k}})^{-1}\cdot W$ for
$\omega\in[w]$, and since $w$ is an arbitrary element of $\Lambda^{k}$.
$\square$$\newline$$\newline$\textbf{\emph{Proof of Proposition
\ref{P3}:}} Recall that $\mathcal{P}=\{[\lambda]\::\:\lambda\in\Lambda\}$.
For $w\in\Lambda^{*}$ set $K_{w}=\varphi_{w}(K)$. Define $I:X\rightarrow\mathbb{R}$
by $I(\omega,W)=-\log\mu_{W,\omega}[\omega_{0}]$ for $(\omega,W)\in X$.
It follows from property (c) stated in Section \ref{S1.2}, from the
ergodic theorem, and from Lemma \ref{L19}, that for $\nu$-a.e. $(\omega,W)\in X$
\begin{multline}
\int H_{\mu}(\mathcal{P}\mid\mathcal{F}_{U})\:d\mu_{F}(U)\\
=\int\int-\log E_{\mu}[1_{[\eta_{0}]}\mid\mathcal{F}_{U}](\eta)\:d\mu(\eta)\:d\mu_{F}(U)\\
=\int\int-\log\mu_{U,\eta}[\eta_{0}]\:d\mu(\eta)\:d\mu_{F}(U)=\int I(\eta,U)\:d\nu(\eta,U)\\
=\underset{n}{\lim}\frac{1}{n}\sum_{k=0}^{n-1}I\circ T^{k}(\omega,W)=\underset{n}{\lim}-\frac{1}{n}\sum_{k=0}^{n-1}\log\mu_{(A_{\omega|_{k}})^{-1}\cdot W,\sigma^{k}\omega}[\omega_{k}]\\
=\underset{n}{\lim}-\frac{1}{n}\sum_{k=0}^{n-1}\log\frac{\mu_{W,\omega}[\omega|_{k+1}]}{\mu_{W,\omega}[\omega|_{k}]}=\underset{n}{\lim}\frac{-\log\mu_{W,\omega}[\omega|_{n}]}{n}\\
=\underset{n}{\lim}\frac{-\log\pi\mu_{W,\omega}(K_{\omega|_{n}})}{n}\:.\label{E6}
\end{multline}
Let $0<\epsilon<-\gamma_{1}$, then there exists a Borel set $\Omega_{0}\in\Omega$
with $\mu(\Omega\setminus\Omega_{0})=0$, such that for $\omega\in\Omega_{0}$
there exists $N_{\omega}\ge1$ for which
\[
\alpha_{i}(A_{\omega|_{n}})\in(e^{n(\gamma_{i}-\epsilon)},e^{n(\gamma_{i}+\epsilon)})\mbox{ for }n\ge N_{\omega}\mbox{ and }1\le i\le d\:.
\]
Since $v\in\mathcal{V}$ there exists $\rho>0$ with 
\[
\rho<\min\{d(\varphi_{\lambda_{1}}(K),\varphi_{\lambda_{2}}(K))\::\:\lambda_{1},\lambda_{2}\in\Lambda\mbox{ with }\lambda_{1}\ne\lambda_{2}\}\:.
\]
Let $\omega\in\Omega_{0}$, $n\ge N_{\omega}$, and $\lambda_{0}\cdot...\cdot\lambda_{n-1}=w\in\Lambda^{n}\setminus\{\omega|_{n}\}$.
Let $0\le k<n$ be such that $\lambda_{k}\ne\omega_{k}$ with $\lambda_{j}=\omega_{j}$
for $0\le j<k$. Since $\pi(\sigma^{k}\omega)\in K_{\omega_{k}}$
we have $B(\pi(\sigma^{k}\omega),\rho)\cap K_{\lambda_{k}}=\emptyset$,
and so
\[
\emptyset=\varphi_{\omega|_{k}}(B(\pi(\sigma^{k}\omega),\rho)\cap K_{\lambda_{k}})\supset\varphi_{\omega|_{k}}(B(\pi(\sigma^{k}\omega),\rho))\cap K_{w}\:.
\]
Now since
\begin{multline*}
\varphi_{\omega|_{k}}(B(\pi(\sigma^{k}\omega),\rho))\supset B(\varphi_{\omega|_{k}}\circ\pi(\sigma^{k}\omega),\alpha_{d}(A_{\omega|_{k}})\cdot\rho)\\
\supset B(\pi(\omega),\alpha_{d}(A_{\omega|_{n}})\cdot\rho)\supset B(\pi(\omega),e^{n(\gamma_{d}-\epsilon)}\cdot\rho),
\end{multline*}
we get $B(\pi(\omega),e^{n(\gamma_{d}-\epsilon)}\cdot\rho)\cap K_{w}=\emptyset$.
We have thus shown that
\[
B(\pi(\omega),e^{n(\gamma_{d}-\epsilon)}\cdot\rho)\cap K_{w}=\emptyset\mbox{ for }\omega\in\Omega_{0},\:n\ge N_{\omega},\mbox{ and }w\in\Lambda^{n}\setminus\{\omega|_{n}\}\:.
\]
It follows from this, from the fact that $\pi\mu_{W,\omega}$ is supported
on $K$ for $\nu$-a.e. $(\omega,W)\in X$, and from (\ref{E6}),
that for $\nu$-a.e. $(\omega,W)\in X$
\begin{multline}
\underset{\delta\downarrow0}{\liminf}\:\frac{\log(\pi\mu_{W,\omega}(B(\pi(\omega),\delta)))}{\log\delta}\\
=\underset{n}{\liminf}\:\frac{\log(\pi\mu_{W,\omega}(B(\pi(\omega),\rho\cdot e^{n(\gamma_{d}-\epsilon)})\cap K))}{\log(\rho\cdot e^{n(\gamma_{d}-\epsilon)})}\\
\ge\underset{n}{\lim}\:\frac{\log(\pi\mu_{W,\omega}(K_{\omega|_{n}}))}{n\cdot(\gamma_{d}-\epsilon)}=\frac{\int H_{\mu}(\mathcal{P}\mid\mathcal{F}_{U})\:d\mu_{F}(U)}{\epsilon-\gamma_{d}}\:.\label{E8}
\end{multline}
For each $\omega\in\Omega$ and $n\ge1$ set $D_{n,\omega}=diag(\alpha_{1}(A_{\omega|_{n}}),...,\alpha_{d}(A_{\omega|_{n}}))$,
let $U_{n,\omega}D_{n,\omega}V_{n,\omega}$ be a singular value decomposition
of $A_{\omega|_{n}}$, and set $L_{n,\omega}=span\{U_{n,\omega}e_{d-m+1},...,U_{n,\omega}e_{d}\}$.
From Lemma \ref{L16} it follows that for $\mu$-a.e. $\omega\in\Omega$
there exists $L_{\omega}\in G_{d,m}$ such that $\{L_{n,\omega}\}_{n=1}^{\infty}$
converges to $L_{\omega}$ in $G_{d,m}$. Set
\[
X_{0}=\{(\omega,W)\in X\::\:\omega\in\Omega_{0},\mbox{ the limit }L_{\omega}=\underset{n}{\lim}\:L_{n,\omega}\mbox{ exists, and }L_{\omega}^{\perp}+W=\mathbb{R}^{d}\},
\]
and for $U\in G_{d,m}$ set
\[
\mathcal{S}_{U}=\{W\in G_{d,m}\::\:U^{\perp}+W\ne\mathbb{R}^{d}\}\:.
\]
From Fubini's theorem and Lemma \ref{L17} we get
\[
\nu(X\setminus X_{0})\le\int_{\{L_{\omega}\mbox{ exists}\}}\mu_{F}(\mathcal{S}_{L_{\omega}})\:d\mu(\omega)=0\:.
\]
Let $b\in(0,\infty)$ be such that $K\subset B(0,b)$. Fix $(\omega,W)\in X_{0}$,
then $L_{\omega}^{\perp}\cap W=\{0\}$, so $P_{L_{\omega}}(x)\ne0$
for each $x\in W\setminus\{0\}$, and so
\[
a_{\omega,W}:=\min\{|P_{L_{\omega}}(x)|\::\:x\in W\mbox{ and }|x|=1\}>0\:.
\]
Since $\{L_{n,\omega}\}_{n=1}^{\infty}$ converges to $L_{\omega}$
it follows that there exists $N_{\omega.W}\ge N_{\omega}$ with
\[
\min\{|P_{L_{n,\omega}}(x)|\::\:x\in W\mbox{ and }|x|=1\}>\frac{a_{\omega,W}}{2}\:\mbox{ for every }n\ge N_{\omega.W}\:.
\]
Let $n\ge N_{\omega.W}$, and set
\[
R=\pi(\omega)+L_{n,\omega}^{\perp}+\{x\in L_{n,\omega}\::\:|x|\le2b\cdot e^{n(\gamma_{d}+\epsilon)}\}\:.
\]
For $d-m+1\le i\le d$ we have $\gamma_{i}=\gamma_{d}$, hence $\alpha_{i}(A_{\omega|_{n}})\le e^{n(\gamma_{d}+\epsilon)}$,
and so
\begin{multline*}
A_{\omega|_{n}}(B(0,2b))=U_{n,\omega}D_{n,\omega}V_{n,\omega}(B(0,2b))=U_{n,\omega}D_{n,\omega}(B(0,2b))\\
\subset U_{n,\omega}(span\{e_{1},...,e_{d-m}\}+\{x\in span\{e_{d-m+1},...,e_{d}\}\::\:|x|\le2b\cdot e^{n(\gamma_{d}+\epsilon)}\})\\
=L_{n,\omega}^{\perp}+\{x\in L_{n,\omega}\::\:|x|\le2b\cdot e^{n(\gamma_{d}+\epsilon)}\}\:.
\end{multline*}
It follows that for $y\in K$
\begin{multline*}
\varphi_{\omega|_{n}}(y)-\pi(\omega)=\varphi_{\omega|_{n}}(y)-\varphi_{\omega|_{n}}\circ\pi\circ\sigma^{n}(\omega)\\
=A_{\omega|_{n}}(y-\pi\circ\sigma^{n}(\omega))\in A_{\omega|_{n}}(B(0,2b))\\
\subset L_{n,\omega}^{\perp}+\{x\in L_{n,\omega}\::\:|x|\le2b\cdot e^{n(\gamma_{d}+\epsilon)}\},
\end{multline*}
which shows that $K_{\omega|_{n}}\subset R$. Given $x\in W$ with
$|x|>\frac{4b}{a_{\omega,W}}\cdot e^{n(\gamma_{d}+\epsilon)}$ we
have
\[
|P_{L_{n,\omega}}(x)|=|x|\cdot|P_{L_{n,\omega}}(\frac{x}{|x|})|>|x|\cdot\frac{a_{\omega,W}}{2}>2b\cdot e^{n(\gamma_{d}+\epsilon)}\:.
\]
It follows that $x+\pi(\omega)\notin R$, and so
\[
(\pi(\omega)+W)\cap K_{\omega|_{n}}\subset(\pi(\omega)+W)\cap R\subset B(\pi(\omega),\frac{4b}{a_{\omega,W}}\cdot e^{n(\gamma_{d}+\epsilon)})\:.
\]
We have thus shown that
\begin{equation}
K_{\omega|_{n}}\cap(\pi(\omega)+W)\subset B(\pi(\omega),\frac{4b}{a_{\omega,W}}\cdot e^{n(\gamma_{d}+\epsilon)})\mbox{ for every }(\omega,W)\in X_{0}\mbox{ and }n\ge N_{\omega,W}\:.\label{E9}
\end{equation}
From property (a) stated in Section \ref{S1.2} it follows that $\pi\mu_{W,\omega}$
is supported on $\pi(\omega)+W$ for $\nu$-a.e. $(\omega,W)\in X$.
From this, from (\ref{E9}), and from (\ref{E6}), we get that for
$\nu$-a.e. $(\omega,W)\in X$
\begin{multline}
\underset{\delta\downarrow0}{\limsup}\:\frac{\log(\pi\mu_{W,\omega}(B(\pi(\omega),\delta)))}{\log\delta}\\
=\underset{n}{\limsup}\:\frac{\log(\pi\mu_{W,\omega}(B(\pi(\omega),\frac{4b}{a_{\omega,W}}\cdot e^{n(\gamma_{d}+\epsilon)})))}{\log(\frac{4b}{a_{\omega,W}}\cdot e^{n(\gamma_{d}+\epsilon)})}\\
\le\underset{n}{\lim}\:\frac{\log(\pi\mu_{W,\omega}(K_{\omega|_{n}}\cap(\pi(\omega)+W)))}{n\cdot(\gamma_{d}+\epsilon)}\\
=\underset{n}{\lim}\:\frac{\log(\pi\mu_{W,\omega}(K_{\omega|_{n}}))}{n\cdot(\gamma_{d}+\epsilon)}=\frac{\int H_{\mu}(\mathcal{P}\mid\mathcal{F}_{U})\:d\mu_{F}(U)}{-\gamma_{d}-\epsilon}\:.\label{E10}
\end{multline}
Now since $\epsilon>0$ can be chosen arbitrarily small the proposition
follows from (\ref{E8}) and (\ref{E10}). $\square$

\section{\label{S6}Proofs of auxiliary Lemmas}

\textbf{\emph{Proof of Lemma \ref{L6}:}} Given a continuous $g:\mathbb{R}^{d}\rightarrow\mathbb{R}$
with compact support it holds for $\mu$-a.e. $\omega$ that
\begin{multline*}
\int g\:d(\pi_{v}\mu)_{W,\pi_{v}(\omega)}=\underset{\delta\downarrow0}{\lim}\frac{1}{P_{W^{\perp}}\pi_{v}\mu(B(P_{W^{\perp}}\pi_{v}(\omega),\delta))}\cdot\int_{P_{W^{\perp}}^{-1}(B(P_{W^{\perp}}\pi_{v}(\omega),\delta))}g\:d\pi_{v}\mu\\
=\underset{\delta\downarrow0}{\lim}\frac{1}{P_{W^{\perp}}\pi_{v}\mu(B(P_{W^{\perp}}\pi_{v}(\omega),\delta))}\cdot\int_{\pi_{v}^{-1}\circ P_{W^{\perp}}^{-1}(B(P_{W^{\perp}}\pi_{v}(\omega),\delta))}g\circ\pi_{v}\:d\mu\\
=\int g\circ\pi_{v}\:d\mu_{v,W,\omega}=\int g\:d\pi_{v}\mu_{v,W,\omega},
\end{multline*}
which proves the Lemma. $\square$

\textbf{\emph{Proof of Lemma \ref{L7}:}} Fix $W\in G_{d,m}$ and
$v_{0}\in\mathcal{V}$, and for each $v\in\mathcal{V}$ set $F_{W}(v)=H_{\mu}(\mathcal{P}\mid\mathcal{F}_{v,W})$,
then it suffice to show that $F_{W}:\mathcal{V}\rightarrow\mathbb{R}$
is upper semi-continuous at $v_{0}$. Let $\{u_{1},...,u_{d-m}\}$
be an orthonormal basis for $W^{\perp}$, and for $1\le i\le d-m$
set $U_{i}=span\{u_{i}\}$ and
\[
\mathcal{Q}_{i}=\{t\in\mathbb{R}\::\:P_{U_{i}}\pi_{v_{0}}\mu\{t\cdot u_{i}\}=0\}\:.
\]
Clearly $\mathbb{R}\setminus\mathcal{Q}_{i}$ is at most countable.
For each $1\le i\le d-m$ and $n\ge1$ let $\{a_{n,k}^{i}\}_{k=-\infty}^{\infty}=\mathcal{J}_{n}^{i}\subset Q_{i}$
be such that $2^{-n-1}\le a_{n,k+1}^{i}-a_{n,k}^{i}\le2^{-n}$ for
$k\in\mathbb{Z}$, and such that $\mathcal{J}_{n}^{i}\subset\mathcal{J}_{n+1}^{i}$.
For $n\ge1$ and $(k_{1},...,k_{d-m})=\bar{k}\in\mathbb{Z}^{d-m}$
set
\[
S_{n,\bar{k}}=P_{W^{\perp}}^{-1}\{\sum_{i=1}^{d-m}t^{i}\cdot u_{i}\::\:(t^{1},...,t^{d-m})\in[a_{n,k_{1}}^{1},a_{n,k_{1}+1}^{1})\times...\times[a_{n,k_{d-m}}^{d-m},a_{n,k_{d-m}+1}^{d-m})\}\:.
\]
For $n\ge1$ and $v\in\mathcal{V}$ let $\mathcal{G}_{v,n}$ be the
$\sigma$-algebra on $\Omega$ generated by 
\[
\{\pi_{v}^{-1}(S_{n,\bar{k}})\::\:\bar{k}\in\mathbb{Z}^{d-m}\},
\]
and set $F_{W,n}(v)=H_{\mu}(\mathcal{P}\mid\mathcal{G}_{v,n})$. For
$v\in\mathcal{V}$ we have $\mathcal{G}_{v,1}\subset\mathcal{G}_{v,2}\subset...$
and $\mathcal{F}_{v,W}=\bigvee_{n=1}^{\infty}\mathcal{G}_{v,n}$,
hence from Theorem 6 in page 38 of \cite{P} we get that $F_{W,1}\ge F_{W,2}\ge...$
and $F_{W}=\underset{n}{\lim}\:F_{W,n}$. It follows that it is enough
to prove that $F_{W,n}:\mathcal{V}\rightarrow\mathbb{R}$ is continuous
at $v_{0}$ for $n\ge1$. Let $n\ge1$, $(k_{1},...,d_{d-m})=\bar{k}\in\mathbb{Z}^{d-m}$
and $\lambda\in\Lambda$ be given, and for $v\in\mathcal{V}$ set
$f(v)=\mu([\lambda]\cap\pi_{v}^{-1}(S_{n,\bar{k}}))$. From the way
$F_{W,n}$ is defined it follows that it suffice to show that $f$
is continuous at $v_{0}$. From $a_{n,k_{i}}^{i},a_{n,k_{i}+1}^{i}\in Q_{i}$
for each $1\le i\le d-m$ it follows that $\mu(\pi_{v_{0}}^{-1}(\partial S_{n,\bar{k}}))=0$,
and for $\omega\in\Omega\setminus\pi_{v_{0}}^{-1}(\partial S_{n,\bar{k}})$
we have
\[
\underset{v\rightarrow v_{0}}{\lim}1_{[\lambda]\cap\pi_{v}^{-1}(S_{n,\bar{k}})}(\omega)=1_{[\lambda]\cap\pi_{v_{0}}^{-1}(S_{n,\bar{k}})}(\omega),
\]
hence from the dominated convergence theorem $\underset{v\rightarrow v_{0}}{\lim}f(v)=f(v_{0})$.
This completes the proof of the lemma. $\square$

\textbf{\emph{Proof of Lemma \ref{L8}:}} Since $\pi_{v}\mu$ is supported
on $K_{v}$ it suffice to show that $\mathcal{L}eb_{d}(K_{v})=0$.
Let $\rho>0$ be such that
\[
\rho<\frac{1}{2}\cdot\min\{d(\varphi_{v,\lambda_{1}}(K_{v}),\varphi_{v,\lambda_{2}}(K_{v}))\::\:\lambda_{1},\lambda_{2}\in\Lambda\mbox{ with }\lambda_{1}\ne\lambda_{2}\}
\]
and set $U=\{x\in\mathbb{R}^{d}\::\:d(x,K_{v})<\rho\}$, then $\varphi_{v,\lambda_{1}}(U)\subset U$
and $\varphi_{v,\lambda_{1}}(U)\cap\varphi_{v,\lambda_{2}}(U)=\emptyset$
for $\lambda_{1},\lambda_{2}\in\Lambda$ with $\lambda_{1}\ne\lambda_{2}$.
Also it is easy to see that the set $U\setminus\cup_{\lambda\in\Lambda}\varphi_{v,\lambda}(U)$
has a non empty interior, hence
\[
\mathcal{L}eb_{d}(U)>\mathcal{L}eb_{d}(\cup_{\lambda\in\Lambda}\varphi_{v,\lambda}(U))=\sum_{\lambda\in\Lambda}\mathcal{L}eb_{d}(\varphi_{v,\lambda}(U))=\mathcal{L}eb_{d}(U)\cdot\sum_{\lambda\in\Lambda}|\det(A_{\lambda})|,
\]
and so $\sum_{\lambda\in\Lambda}|\det(A_{\lambda})|<1$. In addition,
for each $n\ge1$ we have
\begin{multline*}
\mathcal{L}eb_{d}(K_{v})\le\mathcal{L}eb_{d}(\cup_{w\in\Lambda^{n}}\varphi_{v,w}(U))=\sum_{w\in\Lambda^{n}}\mathcal{L}eb_{d}(\varphi_{v,w}(U))\\
=\mathcal{L}eb_{d}(U)\cdot\sum_{w\in\Lambda^{n}}|\det(A_{w})|=\mathcal{L}eb_{d}(U)\cdot\sum_{\lambda_{1},...,\lambda_{n}\in\Lambda}\:\prod_{i=1}^{n}|\det(A_{\lambda_{i}})|\\
=\mathcal{L}eb_{d}(U)\cdot(\sum_{\lambda\in\Lambda}|\det(A_{\lambda})|)^{n},
\end{multline*}
which shows that $\mathcal{L}eb_{d}(K_{v})=0$. $\square$

For the proof of Lemma \ref{L9} we shall first need the following
Lemma regarding the dimension of exceptional sets of projections.
Given $\theta\in\mathcal{M}(\mathbb{R}^{d})$ and $t>0$ let $I_{t}(\theta)$
be the $t$-energy of $\theta$ (see Section 2.5 of \cite{M2}), and
let $\dim_{S}\theta$ be the Sobolev dimension of $\theta$ (see Section
5.2 of \cite{M2}). Given a Borel set $E\subset\mathbb{R}^{d}$ we
denote the restriction of $\theta$ to $E$ by $\theta|_{E}$.
\begin{lem}
\label{L20}Let $\theta\in\mathcal{M}(\mathbb{R}^{d})$ and $1\le l<d$
be given and set $s=\dim_{H}\theta$, then:

(a) If $s\le l$ then for $0<t\le s$
\[
\dim_{H}\{W\in G_{d,l}\::\:\dim_{H}(P_{W}\theta)<t\}\le l(d-l-1)+t\:.
\]

(b) If $s>l$ then for $s-l(d-l)\le t\le l$
\[
\dim_{H}\{W\in G_{d,l}\::\:\dim_{H}(P_{W}\theta)<t\}\le l(d-l)+t-s\:.
\]

(c) If $s>l$ then
\[
\dim_{H}(G_{d,l}\setminus\{W\in G_{d,l}\::\:P_{W}\theta\ll\mathcal{H}^{l}\})\le l(d-l+1)-s,
\]
where $\mathcal{H}^{l}$ is the $l$-dimensional Hausdorff measure.
\end{lem}

\textbf{\emph{Proof of Lemma \ref{L20}, part (a):}} Let $0<t_{0}<t_{1}<t$,
and for each $n\ge1$ set
\[
E_{n}=\{x\in\mathbb{R}^{d}\::\:\theta(B(x,\delta))\le n\cdot\delta^{t_{1}}\mbox{ for each }\delta>0\}\:.
\]
From $\dim_{H}\theta>t_{1}$ and (\ref{E2}) we get $\theta(\mathbb{R}^{d}\setminus\cup_{n}E_{n})=0$.
From an argument as the one given in page 19 of \cite{M2} it follows
that $I_{t_{0}}(\theta|_{E_{n}})<\infty$ for each $n\ge1$. From
this, from Theorem 5.10 in \cite{M2}, and since $\dim_{S}\xi\le\dim_{H}\xi$
for each $\xi\in\mathcal{M}(\mathbb{R}^{d})$ with $\dim_{S}\xi\le d$,
we get
\begin{multline*}
\dim_{H}\{W\in G_{d,l}\::\:\dim_{H}(P_{W}\theta)<t_{0}\}\\
=\underset{n\ge1}{\sup}\:\dim_{H}\{W\in G_{d,l}\::\:\dim_{H}(P_{W}(\theta|_{E_{n}}))<t_{0}\}\le l(d-l-1)+t\;.
\end{multline*}
As this holds for every $0<t_{0}<t$ we obtain \textbf{a} .

\textbf{\emph{Proof of part (b):}} Let $l<t_{0}<t_{1}<s$, and for
each $n\ge1$ let $E_{n}$ be as in the proof of \textbf{a}. Since
$I_{t_{0}}(\theta|_{E_{n}})<\infty$ for each $n\ge1$, it follows
from Theorem 5.10 in \cite{M2} that
\begin{multline*}
\dim_{H}\{W\in G_{d,l}\::\:\dim_{H}(P_{W}\theta)<t\}\\
=\underset{n\ge1}{\sup}\:\dim_{H}\{W\in G_{d,l}\::\:\dim_{H}(P_{W}(\theta|_{E_{n}}))<t\}\le l(d-l)+t-t_{0}\:.
\end{multline*}
Now by letting $t_{0}$ tend to $s$ we obtain \textbf{b}.

\textbf{\emph{Proof of part (c):}} Let $l<t_{2}<t_{0}<t_{1}<s$, and
for each $n\ge1$ let $E_{n}$ be as in the proof of \textbf{a}. Since
$I_{t_{0}}(\theta|_{E_{n}})<\infty$ for each $n\ge1$, it follows
from Theorems 5.4.b and 5.10 in \cite{M2} that
\begin{multline*}
\dim_{H}(G_{d,l}\setminus\{W\in G_{d,l}\::\:P_{W}\theta\ll\mathcal{H}^{l}\})\\
=\underset{n\ge1}{\sup}\:\dim_{H}(G_{d,l}\setminus\{W\in G_{d,l}\::\:P_{W}(\theta|_{E_{n}})\ll\mathcal{H}^{l}\})\\
\le\underset{n\ge1}{\sup}\:\dim_{H}\{W\in G_{d,l}\::\:\dim_{S}(P_{W}(\theta|_{E_{n}}))<t_{2}\}\le l(d-l)+t_{2}-t_{0}\:.
\end{multline*}
Now by letting $t_{2}$ tend to $l$ and $t_{0}$ tend to $s$ we
obtain \textbf{c}. $\square$$\newline$

For the proof of Lemma \ref{L9} we shall also need the following
proposition, which follows directly from Theorem 5.8 in \cite{F2}.
The proof is actually given in \cite{F2} for the case $d=2$, but
extends to higher dimensions without difficulty.
\begin{prop}
\label{P21}Let $1\le l<d$, $E\subset\mathbb{R}^{d}$, $W\in G_{d,l}$,
$\emptyset\ne A\subset W^{\perp}$, and $t>0$ be given. If $\dim_{H}(E\cap(x+W))\ge t$
for each $x\in A$, then $\dim_{H}E\ge t+\dim_{H}A$.
\end{prop}

\textbf{\emph{Proof of Lemma \ref{L9}, part (a):}}\emph{ }Assume
by contradiction that the claim is false for some $0<t\le s$, then
\begin{equation}
\dim_{H}\{W\in G_{d,l}\::\:\mathrm{essinf}_{\theta}\{\dim_{H}(\theta_{W,x})\::\:x\in\mathbb{R}^{d}\}>s-t\}>(l-1)(d-l)+t\:.\label{E17}
\end{equation}
Since the map that sends $W\in G_{d,l}$ to $W^{\perp}\in G_{d,d-l}$
is an isometry with respect to the metric on the Grassmannian defined
in Section \ref{S2}, we get from part (a) of Lemma \ref{L20} that
\begin{multline*}
\dim_{H}\{W\in G_{d,l}\::\:\dim_{H}(P_{W^{\perp}}\theta)<t\}\\
=\dim_{H}\{W\in G_{d,d-l}\::\:\dim_{H}(P_{W}\theta)<t\}\le(l-1)(d-l)+t\:.
\end{multline*}
From this and (\ref{E17}) it follows that there exists $0<\epsilon<t$
and $W\in G_{d,l}$ such that $\dim_{H}(P_{W^{\perp}}\theta)\ge t$
and
\[
\mathrm{essinf}_{\theta}\{\dim_{H}(\theta_{W,x})\::\:x\in\mathbb{R}^{d}\}>s-t+\epsilon\:.
\]
Let $E\subset\mathbb{R}^{d}$ be a Borel set with $\theta(E)>0$,
for $x\in W^{\perp}$ set $E_{x}=E\cap(x+W)$, and set
\[
A=\{x\in W^{\perp}\::\:\theta_{W,x}(E_{x})>0\mbox{ and }\dim_{H}(\theta_{W,x})\ge s-t+\epsilon\}\:.
\]
From properties stated in Section \ref{S1.2} it follows that $P_{W^{\perp}}\theta(A)>0$,
hence 
\[
\dim_{H}A\ge\dim_{H}(P_{W^{\perp}}\theta)\ge t\:.
\]
For $x\in A$ we have
\[
\dim_{H}E_{x}\ge\dim_{H}(\theta_{W,x})\ge s-t+\epsilon,
\]
and so from Proposition \ref{P21} we obtain $\dim_{H}E\ge s+\epsilon$.
As this holds for every Borel set $E\subset\mathbb{R}^{d}$ with $\theta(E)>0$,
it follows that $s=\dim_{H}\theta\ge s+\epsilon$. This is clearly
a contradiction, and so we obtain part (a) of the lemma. The proof
of part (b) is the same, except we need to use part (b) of Lemma \ref{L20}
instead of part (a).\textbf{\emph{$\newline$$\newline$Proof of part
(c):}} Set
\[
S=\{W\in G_{d,l}\::\:P_{W^{\perp}}\theta\ll\mathcal{H}^{d-l}\},
\]
then from part (c) of Lemma \ref{L20} we get 
\begin{equation}
\dim_{H}(G_{d,l}\setminus S)\le(d-l)(l+1)-s\:.\label{E18}
\end{equation}
Let $d-l<t_{0}<t_{1}<s$ and for $n\ge1$ set
\[
E_{n}=\{x\in\mathbb{R}^{d}\::\:\theta(B(x,\delta))\le n\cdot\delta^{t_{1}}\mbox{ for each }\delta>0\},
\]
then as in the proof of part (a) of Lemma \ref{L20} we have $\theta(\mathbb{R}^{d}\setminus\cup_{n}E_{n})=0$
and $I_{t_{0}}(\theta|_{E_{n}})<\infty$ for each $n\ge1$. Since
for each $W\in G_{d,l}$ we have $\theta_{W,x}(\mathbb{R}^{d}\setminus\cup_{n}E_{n})=0$
for $\theta$-a.e. $x\in\mathbb{R}^{d}$, it follows that
\begin{multline}
\dim_{H}\{W\in S\::\:\mathrm{essinf}_{\theta}\{\dim_{H}(\theta_{W,x})\::\:x\in\mathbb{R}^{d}\}<t_{0}-d+l\}\\
=\underset{n\ge1}{\sup}\:\dim_{H}\{W\in S\::\:\mathrm{essinf}_{\theta}\{\dim_{H}(\theta_{W,x}|_{E_{n}})\::\:x\in\mathbb{R}^{d}\}<t_{0}-d+l\}\:.\label{E5}
\end{multline}
As described in Section 2 of \cite{JM}, given $W\in G_{d,l}$ and
a Radon measure $\xi$ on $\mathbb{R}^{d}$ with compact support,
there exist Radon measures $\{\xi^{W,x}\}_{x\in W^{\perp}}$ on $\mathbb{R}^{d}$
such that for $\mathcal{H}^{d-l}$-a.e. $x\in W^{\perp}$
\[
\int g\:d\xi^{W,x}=\underset{\delta\downarrow0}{\lim}\:\frac{1}{(2\delta)^{d-l}}\cdot\int_{P_{W^{\perp}}^{-1}(B(x,\delta))}g\:d\xi\;\mbox{ for }g\in C(\mathbb{R}^{d})\:.
\]
For $x\in\mathbb{R}^{d}$ we set $\xi^{W,x}:=\xi^{W,P_{W^{\perp}}x}$.$\newline$Fix
some $n\ge1$ with $\theta(E_{n})>0$, and let $W\in S$. From property
(b) in Section \ref{S1.2} above and from Theorem 2.12 in \cite{M3},
it follows that for $\theta$-a.e. $x\in\mathbb{R}^{d}$ we have for
each $g\in C(\mathbb{R}^{d})$
\begin{multline*}
\int g\:d\theta^{W,x}=\underset{\delta\downarrow0}{\lim}\:\frac{P_{W^{\perp}}\theta(B(P_{W^{\perp}}x,\delta))}{(2\delta)^{d-l}}\cdot\frac{\int_{P_{W^{\perp}}^{-1}(B(P_{W^{\perp}}x,\delta))}g\:d\theta}{P_{W^{\perp}}\theta(B(P_{W^{\perp}}x,\delta))}\\
=\frac{dP_{W^{\perp}}\theta}{d\mathcal{H}^{d-l}}(P_{W^{\perp}}x)\cdot\int g\:d\theta_{W,x},
\end{multline*}
which shows that
\[
\theta^{W,x}=\frac{dP_{W^{\perp}}\theta}{d\mathcal{H}^{d-l}}(P_{W^{\perp}}x)\cdot\theta_{W,x}\:.
\]
From this, from $0<\frac{dP_{W^{\perp}}\theta}{d\mathcal{H}^{d-l}}(P_{W^{\perp}}x)<\infty$
for $\theta$-a.e. $x\in\mathbb{R}^{d}$, and from Lemma 3.2 in \cite{JM},
we get that for $\theta$-a.e. $x\in\mathbb{R}^{d}$
\[
\dim_{H}(\theta_{W,x}|_{E_{n}})=\dim_{H}(\theta^{W,x}|_{E_{n}})=\dim_{H}((\theta|_{E_{n}})^{W,x})\:.
\]
Now from Lemma 2.22 in \cite{JM}, from $I_{t_{0}}(\theta|_{E_{n}})<\infty$,
and from Theorem 6.5 in \cite{M2}, we obtain
\begin{multline*}
\dim_{H}\{W\in S\::\:\mathrm{essinf}_{\theta}\{\dim_{H}(\theta_{W,x}|_{E_{n}})\::\:x\in\mathbb{R}^{d}\}<t_{0}-d+l\}\\
=\dim_{H}\{W\in S\::\:\mathrm{essinf}_{\theta}\{\dim_{H}((\theta|_{E_{n}})^{W,x})\::\:x\in\mathbb{R}^{d}\}<t_{0}-d+l\}\\
\le\dim_{H}\{W\in S\::\:\int_{W^{\perp}}I_{t_{0}-d+l}((\theta|_{E_{n}})^{W,x})\:d\mathcal{H}^{d-l}(x)=\infty\}\le(d-l)(l+1)-t_{0}\:.
\end{multline*}
This together with (\ref{E18}) and (\ref{E5}) proves part (c) of
the lemma, since we can let $t_{0}$ tend to $s$. $\square$$\newline$$\newline$\textbf{\emph{Proof
of Lemma \ref{L10}:}} Fix $v\in\mathbb{R}^{d|\Lambda|}$ and set
$\pi=\pi_{v}$, $K=K_{v}$, and $\varphi_{\lambda}=\varphi_{v,\lambda}$
for $\lambda\in\Lambda$. Let $k:=k(\mu)\ge0$ be as defined in (\ref{E19}).
If $D(\mu)\ge d$ then there is nothing to prove (see Proposition
10.3 in \cite{F}), hence we can assume $D(\mu)<d$, and so $k<d$.
For $1\le i\le k$ and $w\in\Lambda^{*}$ set $d_{i,w}=\left\lceil \frac{\alpha_{i}(A_{w})}{\alpha_{k+1}(A_{w})}\right\rceil $,
and set 
\[
d_{w}=\begin{cases}
\prod_{i=1}^{k}d_{i,w} & \mbox{, if }k>0\\
1 & \mbox{, if }k=0
\end{cases}\:.
\]
There exists a constant $a>0$ such that for each $w\in\Lambda^{*}$
there exists a rectangle $R_{w}\subset\mathbb{R}^{d}$ with $\varphi_{w}(K)\subset R_{w}$,
and with side lengths $s_{1},...,s_{d}>0$ where
\[
s_{i}=\begin{cases}
a\cdot\alpha_{k+1}(A_{w})\cdot d_{i,w} & ,\mbox{ if }1\le i\le k\\
a\cdot\alpha_{k+1}(A_{w}) & ,\mbox{ if }k+1\le i\le d
\end{cases}\:.
\]
For $w\in\Lambda^{*}$ let $\mathcal{R}_{w}=\{R_{w,1},...,R_{w,d_{w}}\}$
be a partition of $R_{w}$ into disjoint squares of side length $a\cdot\alpha_{k+1}(A_{w})$.
For $\omega\in\Omega$ and $n\ge1$ let $R_{\omega,n}$ be the unique
member of $\mathcal{R}_{\omega|_{n}}$ which contains $\pi(\omega)$.
For each $n\ge1$ set
\[
E_{n}=\{\omega\in\Omega\::\:\pi\mu(R_{\omega,n})\le\frac{\mu[\omega|_{n}]}{d_{\omega|_{n}}\cdot n^{2}}\},
\]
then
\[
\mu(E_{n})\le\sum_{w\in\Lambda^{n}}\sum_{j=1}^{d_{w}}\pi\mu(R_{w,j})\cdot1_{\{\pi\mu(R_{w,j})\le\frac{\mu[w]}{d_{w}\cdot n^{2}}\}}\le\frac{1}{n^{2}},
\]
and so $\sum_{n=1}^{\infty}\mu(E_{n})<\infty$. From this and the
Borel-Cantelli Lemma it follows that 
\begin{equation}
\mu\{\omega\::\:\#\{n\ge1\::\:\omega\in E_{n}\}=\infty\}=0\:.\label{E20}
\end{equation}
There exists a constant $a'>a$ such that
\[
R_{\omega,n}\subset B(\pi(\omega),a'\cdot\alpha_{k+1}(A_{\omega|_{n}}))\:\mbox{ for }\omega\in\Omega\mbox{ and }n\ge1,
\]
hence for $\omega\in\Omega$
\begin{multline*}
\underset{\delta\downarrow0}{\limsup}\:\frac{\log\pi\mu(B(\pi(\omega),\delta))}{\log\delta}\\
=\underset{n\rightarrow\infty}{\limsup}\:\frac{\log\pi\mu(B(\pi(\omega),a'\cdot\alpha_{k+1}(A_{\omega|_{n}})))}{\log(a'\cdot\alpha_{k+1}(A_{\omega|_{n}}))}\\
\le\underset{n\rightarrow\infty}{\limsup}\:\frac{\log\pi\mu(R_{\omega,n})}{\log(\alpha_{k+1}(A_{\omega|_{n}}))}\:.
\end{multline*}
Now from (\ref{E20}) it follows that for $\mu$-a.e. $\omega\in\Omega$
\begin{multline*}
\underset{\delta\downarrow0}{\limsup}\:\frac{\log\pi\mu(B(\pi(\omega),\delta))}{\log\delta}\le\underset{n\rightarrow\infty}{\limsup}\:\frac{\log(\frac{\mu[\omega|_{n}]}{d_{\omega|_{n}}\cdot n^{2}})}{\log(\alpha_{k+1}(A_{\omega|_{n}}))}\\
=\underset{n\rightarrow\infty}{\limsup}\:\frac{\log\mu[\omega|_{n}]-\sum_{i=1}^{k}\log\frac{\alpha_{i}(A_{\omega|_{n}})}{\alpha_{k+1}(A_{\omega|_{n}})}}{\log(\alpha_{k+1}(A_{\omega|_{n}}))}\:.
\end{multline*}
This together with (\ref{E1}) and the Shannon-McMillan-Breiman theorem
gives
\[
\underset{\delta\downarrow0}{\limsup}\:\frac{\log\pi\mu(B(\pi(\omega),\delta))}{\log\delta}\le k-\frac{h_{\mu}+\gamma_{1}+...+\gamma_{k}}{\gamma_{k+1}}=D(\mu)
\]
for $\mu$-a.e. $\omega\in\Omega$, which proves the lemma. $\square$

\textbf{\emph{Proof of Lemma \ref{L11}:}} Assume by contradiction
that $\mathbf{G}'$ is not $m$-irreducible, then there exists a proper
linear subspace $W$ of $\mathcal{A}^{m}(\mathbb{R}^{d})$ such that
$\mathcal{A}^{m}M(W)=W$ for all $M\in\mathbf{G}'$. Let $W_{1},...,W_{k}$
be an enumeration of the set 
\[
\{\mathcal{A}^{m}A_{w}(W)\::\:w\in\Lambda^{n-1}\}
\]
and define
\[
\mathbf{H}=\{M\in Gl(d,\mathbb{R})\::\;\forall\;1\le i\le k\quad\exists\;1\le j\le k\;\mbox{ with }\mathcal{A}^{m}M(W_{i})=W_{j}\},
\]
then $\mathbf{H}$ is a closed subgroup of $Gl(d,\mathbb{R})$. Let
$\mathbf{T}$ denote the closure of the semigroup generated by $\{A_{\lambda}^{-1}\}_{\lambda\in\Lambda}$.
Since $\mathcal{A}^{m}M(W)=W$ for each $M\in\mathbf{G}'$ it follows
that $\mathbf{H}$ contains the semigroup generated by $\{A_{\lambda}\}_{\lambda\in\Lambda}$,
and so $\mathbf{T}\subset\mathbf{H}$. This implies that $\mathbf{T}$
is not $m$-strongly irreducible which contradicts Lemma \ref{L15},
and so it must hold that $\mathbf{G}'$ is $m$-irreducible.

From Proposition III.5.6 in \cite{BL} it follows that for each $1\le i\le d$
\begin{multline*}
\gamma_{i}'=\underset{N}{\lim}\:\frac{1}{N}\int_{(\Lambda^{n})^{\mathbb{N}}}\log\alpha_{i}(A_{\omega|_{N}})\:d\mu'(\omega)\\
=\underset{N}{\lim\:}\frac{1}{N}\int_{\Lambda^{\mathbb{N}}}\log\alpha_{i}(A_{\omega|_{n\cdot N}})\:d\mu(\omega)=n\cdot\gamma_{i},
\end{multline*}
hence
\[
\max\{1\le i\le d\::\:\gamma_{d-i+1}'=...=\gamma_{d}'\}=m<d\:.
\]
From this, from the $m$-irreducibility of $\mathbf{G}'$, and from
Proposition \ref{L2}, it follows that there exists a unique $\mu_{F}'\in\mathcal{M}(G_{d,m})$
with $\mu_{F}'=\sum_{w\in\Lambda^{n}}p_{w}\cdot A_{w}^{-1}\mu_{F}'$.
Clearly we also have $\mu_{F}=\sum_{w\in\Lambda^{n}}p_{w}\cdot A_{w}^{-1}\mu_{F}$,
hence $\mu_{F}'=\mu_{F}$. $\square$

\end{document}